\documentclass[11pt]{amsart}

\usepackage{amsmath} 
\usepackage{amssymb}

\newtheorem{theorem}{Theorem}[section] 
\newtheorem{claim}{Claim}[theorem]
\newtheorem{mclaim}[claim]{Main Claim}
\newtheorem{proposition}[theorem]{Proposition} 

\theoremstyle{definition}
\newtheorem{definition}[theorem]{Definition}

\theoremstyle{remark}
\newtheorem{remark}[theorem]{Remark}
\newtheorem{context}[theorem]{Context}
\newtheorem{tcontext}[theorem]{Temporary Context}
\newtheorem{notation}[theorem]{Notation}

\newtheorem{observation}[theorem]{Observation}

\newcommand{\bbD}{{\mathbb D}}
\newcommand{\bF}{{\bf F}}
\newcommand{\bbP}{{\mathbb P}}
\newcommand{\bbQ}{{\mathbb Q}}
\newcommand{\bV}{{\bf V}}
\newcommand{\forces}{\Vdash} 
\newcommand{\Dom}{{\rm Dom}}
\newcommand{\Rang}{{\rm Rang}}
\newcommand{\pr}{{\rm pr}}
\newcommand{\st}{{\name{\rm st}}}

\newcommand{\cf}{{\rm cf}}
\newcommand{\cS}{{\mathcal S}}
\newcommand{\cI}{{\mathcal I}}

\newcommand{\otp}{{\rm otp}}
\newcommand{\mrot}{{\rm root}}
\newcommand{\suc}{{\rm succ}}
\newcommand{\splt}{{\rm split}}

\newcount\skewfactor
\def\mathunderaccent#1#2 {\let\theaccent#1\skewfactor#2
\mathpalette\putaccentunder}
\def\putaccentunder#1#2{\oalign{$#1#2$\crcr\hidewidth
\vbox to.2ex{\hbox{$#1\skew\skewfactor\theaccent{}$}\vss}\hidewidth}}
\def\name{\mathunderaccent\tilde-3 }

\title[Iteration of $\lambda$-complete forcing notions]{Iteration of
$\lambda$-complete forcing notions not collapsing $\lambda^+$.}

\author{Andrzej Ros{\l}anowski}
\address{Department of Mathematics\\
 University of Nebraska at Omaha\\
 Omaha, NE 68182-0243, USA} 
\email{roslanowski@unomaha.edu}
\urladdr{http://www.unomaha.edu/$\sim$aroslano}

\author{Saharon Shelah}
\address{Institute of Mathematics\\
 The Hebrew University of Jerusalem\\
 Jerusalem 91904, Israel\\
 and  Department of Mathematics\\
 Rutgers University\\
 New Brunswick, NJ 08854, USA}
\email{shelah@math.huji.ac.il}
\urladdr{http://www.math.rutgers.edu/$\sim$shelah}

\thanks{Both authors would like to thank the United States-Israel
Binational Science Foundation for partial support. This publication has
number 655 in Shelah's list.}


\begin{document}

\begin{abstract}
We look for a parallel to the notion of ``proper forcing'' among
$\lambda$-complete forcing notions not collapsing $\lambda^+$.  We suggest
such a definition and prove that it is preserved by suitable iterations.
\end{abstract}

\maketitle

\setcounter{section}{-1}
\section{Introduction} 
This work follows \cite{Sh:587} and \cite{Sh:667} (and see history there),
but we do not rely on those papers. Our goal in this and the previous papers
is to develop a theory parallel to ``properness in CS iterations'' for
iterations with larger supports. In \cite{Sh:587}, \cite{Sh:667} we have
presented parallels to \cite{Sh:64} and \cite{Sh:98}, whereas here we try to
have parallels to \cite{Sh:100}, \cite[Ch.III]{Sh:b},
\cite[Ch.V,\S5-\S7]{Sh:b} and hopefully \cite[Ch.VI]{Sh:f},
\cite[Ch.XVIII]{Sh:f}. 

It seems too much to hope for a notion fully parallel to ``proper'' among
$\lambda$-complete forcing notions as even for ``$\lambda^+$-c.c.
$\lambda$-complete'' there are problems. We should also remember about ZFC
limitations for possible iteration theorems. For example, if in the
definition of the forcing notion $\bbQ^*$ in Section 3 we demand
$h^p\restriction e_\delta\subseteq h_\delta$, then the proof fails. This may
seem a drawback, but one should look at \cite[AP, p.985, 3.6(2) and p.990,
3.9]{Sh:f}. By it, if $\cS^*=\cS^{\lambda^+}_\lambda$, and ($A_\delta,
h_\delta$ are as in \ref{conlast} and) we ask a success on a club, then for
some $\langle h_\delta:\delta\in \cS^{\lambda^+}_\lambda\rangle$ we fail.
Now, if we allow only $h_\delta:A_\delta\longrightarrow 2$ and we ask for
``success of the uniformization'' on an end segment of $A_\delta$ (for all
such $\langle A_\delta:\delta\in\cS^{\lambda^+}_\lambda\rangle$), then we
also fail as we may code colourings with values in $\lambda$.

In the first section we formulate our definitions (including {\em properness
over $\lambda$}, see \ref{1.2}). We believe that our main Definition
\ref{1.2} is quite reasonable and applicable. One may also define a version
of it where the diamond is ``spread out''. The second section is devoted to
the proof of the preservation theorem, and the next one gives three
(relatively easy) examples of forcing notions fitting our scheme. We
conclude the paper with the discussion of applications and variants.

\begin{notation}
Our notation is rather standard and compatible with that of classical
textbooks (like Jech \cite{J}). In forcing we keep the older (Cohen's)
convention that {\em a stronger condition is the larger one}.    
\begin{enumerate}
\item For a filter $D$ on $\lambda$, the family of all $D$--positive subsets
of $\lambda$ is called $D^+$. (So $A\in D^+$ if and only if $A\subseteq
\lambda$ and $A\cap B\neq\emptyset$ for all $B\in D$.)

\item Every forcing notion $\mathbb P$ under considerations is assumed to
have the weakest condition $\emptyset_{\mathbb P}$, i.e., $(\forall p\in
{\mathbb P})(\emptyset_{\mathbb P}\leq_{\mathbb P} p)$. Also we assume $*
\notin{\mathbb P}$ is a fixed object belonging to all the $N$'s we consider. 

\item A tilde indicates that we are dealing with a name for an object in
a forcing extension (like $\name{x}$). The canonical $\bbP$--name for the
$\bbP$--generic filter over $\bV$ is denoted by $\name{G}_{\bbP}$. In
iterations, if $\bar{\bbQ}=\langle\bbP_\zeta,\name{\bbQ}_\zeta:\zeta<
\zeta^*\rangle$ and $p\in\lim(\bar{\bbQ})$, then we keep convention that
$p(\alpha)=\name{\emptyset}_{\name{\bbQ}_\alpha}$ for $\alpha\in\zeta^*
\setminus\Dom(p)$. 

\item Ordinal numbers will be denoted be the lower case initial letters of
the Greek alphabet ($\alpha,\beta,\gamma\ldots$) and also by $i,j$ (with
possible sub- and superscripts). 

\item A bar above a letter denotes that the object considered is a sequence;
usually $\bar{X}$ will be $\langle X_i:i<\zeta)\rangle$, where $\zeta$
denotes the length of $\bar{X}$. Often our sequences will be indexed by a
set of ordinals, say $\cS\subseteq\lambda$, and then $\bar{X}$ will
typically be $\langle X_\delta:\delta\in\cS\rangle$. Semi--diamond sequences
will be called $\bar{F}$ (with possible superscripts).   
\end{enumerate}
\end{notation}

In our definitions (and proofs) we will use somewhat special diamond--like
sequences (see \ref{pre1.2}(2)). The difference between them and classical
diamonds is quite minor, so let us remind the following.

\begin{definition}
\label{diam}
\begin{enumerate}
\item Let $D$ be a filter on $\lambda$. We say that $\bar{F}=\langle
F_\delta:\delta\in{\mathcal S}\rangle$ is a {\em $D$--diamond sequence} if 
${\mathcal S}\in D^+$, $F_\delta\in {}^\delta\delta$ for $\delta\in
{\mathcal S}$, and 
\[(\forall f\in {}^\lambda \lambda)(\{\delta\in {\mathcal S}:F_\delta
\subseteq f\}\in D^+).\] 
We may also call such $\bar{F}$ {\em a $(D,\cS)$--diamond sequence}.
\item We say that {\em $(D,\cS)$ has diamonds\/} if there is a
$(D,\cS)$--diamond. We say that {\em $D$ has diamonds\/} if $D$ is a normal
filter on $\lambda$ and for every $\cS\in D^+$ there is a
$(D,\cS)$--diamond.
\end{enumerate}
\end{definition}

\begin{definition}
\label{comp}
A forcing notion $\bbP$ is {\em $\lambda$--complete} if every
$\leq_{\bbP}$--increasing chain of length less than $\lambda$ has an upper
bound in $\bbP$. It is {\em $\lambda$--lub--complete} if every
$\leq_{\bbP}$--increasing chain of length less than $\lambda$ has a least
upper bound in $\bbP$.
\end{definition}

\begin{proposition}
\label{basics}
\begin{enumerate}
\item If $D$ is a filter on $\lambda$, then the family of all diagonal
intersections of members of $D$ constitutes a normal filter (but in general 
not necessarily proper). We call this family {\em the normal filter
generated by $D$}. 
\item If ${\mathbb P}$ is a $\lambda$--complete forcing notion and $D$ is a
normal filter on $\lambda$, \underline{then} in $\bV^{\bbP}$ the filter $D$
generates a proper normal filter on $\lambda$. [Abusing notation, we will
denote this filter also by $D$ or, if we want to stress that we work in the
forcing extension, by $D^{\bV[G_\bbP]}$.]  

\noindent Moreover, by the $\lambda$--completeness of $\bbP$, if $X\in
D^+\cap\bV$, then $\forces_{\bbP} X\in D^+$, and if $X\in\bV$,
$p\forces_{\bbP} X\in D^{\bV^\bbP}$ then $X\in D$.  
\item If ${\mathbb P}$ is a $\lambda$--complete forcing notion and $\bar{F}=
\langle F_\delta: \delta\in{\mathcal S}\rangle$ is a $D$--diamond sequence,
then  
\[\forces_{\bbP}\mbox{`` }\bar{F}\mbox{ is a $D$--diamond sequence ''.}\]
\end{enumerate}
\end{proposition}

Definition \ref{pairing} and Proposition \ref{1.1B} below are not central
for us, but they may be used to get somewhat stronger results, see
\cite{Sh:F509}. 

\begin{definition}
\label{pairing}
Let $\pr$ be a definable pairing function on $\lambda$, for example 
$\pr(\alpha,\beta)=\omega^{\alpha+\beta}+\beta$, and let $\bar{F}=\langle
F_\delta:\delta\in {\mathcal S}\rangle$ be a $D$--diamond sequence.

For an ordinal $\alpha<\lambda$ we let $\bar{F}^{[\alpha]}=\langle
F^{[\alpha]}_\delta:\delta\in {\mathcal S}\rangle$, where each
$F^{[\alpha]}_\delta$ is a function with domain $\delta$ and such that 
\[F^{[\alpha]}_\delta(\beta)=\left\{\begin{array}{ll}
F_\delta(\pr(\alpha,\beta))&\mbox{ if well defined},\\
0&\mbox{ otherwise.}
				    \end{array}\right.\]
\end{definition}

\begin{proposition}
\label{1.1B} 
If $\bar{F}$ is a $D$--diamond sequence, \underline{then} for every
$\alpha<\lambda$, $\bar{F}^{[\alpha]}$ is also a $D$--diamond sequence.
\end{proposition}

Throughout the paper we will assume the following:
\begin{context}
\label{1.1} 
\begin{enumerate}
\item[(a)] $\lambda$ is an uncountable cardinal, $\lambda=\lambda^{<
\lambda}$, and 
\item[(b)] $D$ is a normal filter on $\lambda$ (usually $D$ is the club
filter ${\mathcal D}_\lambda$ on $\lambda$), 
\item[(c)] ${\mathcal S}\in D^+$ contains all successor ordinals below
$\lambda$, $0\notin\cS$, and $\cS'=\lambda\setminus\cS$ is unbounded in 
$\lambda$, 
\item[(d)] there is a $(D,\cS)$--diamond sequence.
\end{enumerate}
\end{context}

\section{The definitions} 
In this section we define a special genericity game, properness over $(D,
\cS)$--semi diamonds and the class of forcing notions we are interested in.  

\begin{definition}
\label{pre1.2}
Let $\bbP$ be a forcing notion and let $N\prec ({\mathcal H}(\chi),{\in},
{<^*_\chi})$ be such that $\|N\|=\lambda$, $N^{<\lambda}\subseteq N$ and
$\{\lambda,\mathbb P,D,\cS\}\in N$. Let $h:\lambda\longrightarrow N$ be
such that the range $\Rang(h)$ of the function $h$ includes $\bbP\cap N$.
\begin{enumerate}
\item We say that $\bar{F}=\langle F_\delta:\delta\in {\mathcal S}\rangle$
is {\em a $(D,\cS)$--semi diamond sequence\/} if $F_\delta\in {}^\delta
\delta$ for $\delta\in\cS$ and 
\begin{enumerate}
\item[$(*)$] for every $\leq_\bbP$--increasing sequence $\bar{p}=\langle
p_\alpha:\alpha<\lambda\rangle\subseteq \bbP\cap N$ we have
\[\{\delta\in\cS: (\forall\alpha<\delta)(h\circ F_\delta(\alpha)=p_\alpha)\} 
\in D^+.\]
\end{enumerate}
\item Let $\bar{F}$ be a $(D,\cS)$--semi diamond. A sequence $\bar{q}=
\langle q_\delta:\delta\in\cS\rangle\subseteq N\cap\bbP$ is called {\em an
$(N,h,\bbP)$--candidate over $\bar{F}$} (or: {\em
$(N,h,\bbP,\bar{F})$--candidate}) whenever 
\begin{enumerate}
\item[$(\alpha)$] for every open dense subset ${\mathcal I}\in N$ of $\bbP$
\[\{\delta \in {\mathcal S}:q_\delta\in {\mathcal I}\}={\mathcal S}\mod D,\]
and
\item[$(\beta)$] \underline{if} $\delta\in {\mathcal S}$ is a limit
ordinal and $\langle h\circ F_\delta(\alpha):\alpha<\delta\rangle$ is a
$\le_{\bbP}$--increasing sequence of members of $\bbP\cap N$,\\
\underline{then} $q_\delta$ is its upper bound in $\bbP$.
\end{enumerate}
\item Let $\bar{q}$ be an $(N,h,\bbP,\bar{F})$--candidate and
$r\in\bbP$. We define a game $\Game(r,N,h,\bbP,\bar{F},\bar{q})$ of two
players, the {\em generic player} and the {\em anti-generic player}, as
follows. A play lasts $\lambda$ moves, in the $i^{\rm th}$ move conditions
$r_i^-,r_i\in\bbP$ and a set $C_i\in D$ are chosen such that 
\begin{itemize}
\item $r_i^-\in N$, $r_i^-\leq r_i$, $r\leq r_i$,
\item $(\forall j<i)(r_j\leq r_i\ \&\ r_j^-\leq r_i^-)$, and 
\item the generic player chooses $r_i^-,r_i,C_i$ if $i\in\cS$, and the
anti-generic player chooses $r_i^-,r_i,C_i$ if $i\in {\mathcal S}'$.   
\end{itemize}
If at some moment during the play there is no legal move for one of the
players, then the anti-generic player wins. If the play lasted $\lambda$
moves, then the generic player wins the play whenever
\begin{enumerate}
\item[$(\circledast)$] \underline{if} $\delta\in {\mathcal S}\cap
\bigcap\limits_{i<\delta} C_i$ is a limit ordinal, and $\langle h\circ
F_\delta(\alpha):\alpha<\delta\rangle=\langle r^-_\alpha:\alpha<\delta
\rangle$, \underline{then} $q_\delta\leq r_\delta$.  
\end{enumerate}
\item Let $\bar{q}$ be an $(N,h,\bbP,\bar{F})$--candidate, $\bar{F}$ a
$(D,\cS)$--semi diamond. A condition $r\in{\mathbb P}$ is {\em $(N,h,
\bbP)$--generic for $\bar{q}$ over $\bar{F}$} if the generic player has a
winning strategy in the game $\Game(r,N,h,\bbP,\bar{F},\bar{q})$.
\end{enumerate}
\end{definition}

\begin{observation}
\label{obs}
\begin{enumerate}
\item In the game $\Game(r,N,h,\bbP,\bar{F},\bar{q})$, for each of the
players, if it increases conditions $r_i^-,r_i$, its choice can only improve
its situation. Making sets $C_i$ (for $i\in\cS$) smaller can only help the
generic player.
\item If forcing with ${\mathbb P}$ does not add new subsets to $\lambda$,
then the game in Definition \ref{pre1.2}(5) degenerates as without loss of 
generality $r$ forces a value to $\name{G}_{\mathbb P}\cap N$; the condition
does not degenerate, in fact this condition (which implies adding no new
$\lambda$--sequences) is preserved by $(<\lambda^+)$--support iterations
(see \cite{Sh:587}). 
\item Also if $\cS_1\subseteq\cS\mod D$, $\cS_1\in D^+$, then in Definition
\ref{pre1.2} we can replace $\cS$ by $\cS_1$. (Again, the generic player can
guarantee $C_i\cap\cS_1\subseteq\cS$.) 
\item If $\bbP$ is $\lambda$--complete and $r$ is $(N,\bbP)$--generic (in
the usual sense, i.e., $r\forces_\bbP\mbox{``}\,N[\name{G}_\bbP]\cap\bV=N\,
\mbox{''}$), then both players have always legal moves in the game
$\Game(r,N,h,\bbP,\bar{F},\bar{q})$.\\ 
Also if the forcing notion $\bbP$ is $\lambda$--lub--complete, then both
players have always legal moves in the game $\Game(r,N,h,\bbP,\bar{F},
\bar{q})$ (for any $r$).
\end{enumerate}
\end{observation}

\begin{definition}
\label{1.2} 
\begin{enumerate}
\item Let $\cS\in D^+$. We say that a forcing notion $\bbP$ is {\em proper
over $(D,\cS)$--semi diamonds} whenever (there is a $(D,\cS)$--diamond and): 
\begin{enumerate}
\item[(a)] $\bbP$ is $\lambda$--complete, and 
\item[(b)] \underline{if} $\chi$ is large enough, $p\in\bbP$ and $N\prec
({\mathcal H} (\chi),\in,<^*_\chi)$, $\|N\|=\lambda$, $N^{<\lambda}\subseteq
N$ and $\{\lambda,p,\bbP,D,\cS,\ldots\}\in N$, and $h:\lambda\longrightarrow
N$ satisfies $\bbP\cap N\subseteq\Rang(h)$, and $\bar{F}$ is a
$(D,\cS)$--semi diamond for $(N,h,\bbP)$, and $\bar{q}=\langle q_\delta:
\delta\in {\mathcal S}\rangle$ is an $(N,h,\bbP,\bar{F})$--candidate,  

\underline{then} there is $r\in\bbP$ stronger than $p$ and such that $r$ is
$(N,h,\bbP)$--generic for $\bar{q}$ over $\bar{F}$.
\end{enumerate}
\item $\bbP$ is said to be {\em proper over $D$--semi diamonds} if it is
proper over $(D,{\mathcal S})$--semi diamonds for every ${\mathcal S}\in
D^+$ (so $D$ has diamonds). The family of forcing notions proper over
$D$--semi diamonds is denoted $K^1_D$. 
\item A forcing notion $\bbP$ is {\em proper over $\lambda$} if it is
proper over $D$--semi diamonds for every normal filter $D$ on $\lambda$
which has diamonds.
\end{enumerate}
\end{definition}

\begin{remark}
Does $D$ matter? Yes, as we may use some ``large $D$'' and be interested in
preserving its largeness properties.
\end{remark}

\begin{proposition}
If $\bbP$ is a $\lambda^+$--complete forcing notion, \underline{then}
$\bbP$ is proper over $\lambda$.
\end{proposition}

\begin{proof}
Straightforward.
\end{proof}
 
\begin{proposition}
\label{1.2A} 
\begin{enumerate}
\item If $N,\bbP,h$ are as in \ref{pre1.2}, $\bbP$ is $\lambda$-complete,
and $\bar{F}$ is a $(D,\cS)$--semi diamond, \underline{then} there is an
$(N,h,\bbP,\bar{F})$--candidate. In fact we can even demand: 
\begin{enumerate}
\item[$(+)$] if ${\mathcal I} \in N$ is an open dense subset of $\bbP$,
\underline{then} $q_\delta\in {\mathcal I}$ for every large enough $\delta$. 
\end{enumerate}
\item Let $r$ be $(N,h,\bbP)$--generic over $\bar{F}$ for some $(N,h,\bbP,
\bar{F})$--candidate $\bar{q}$. \underline{Then} 
\begin{enumerate}
\item[(a)] if $\langle r_i^-,r_i,C_i:i<\lambda\rangle$ is a result of a play  
of $\Game(r,N,h,\bbP,\bar{F},\bar{q})$ in which the generic player uses its
winning strategy, then 
\[G'=\{p\in\bbP\cap N:(\exists i<\lambda)(p\leq r_i)\}\]
is a subset of $\bbP\cap N$ generic over $N$, and 
\item[(b)] $r$ is $(N,\bbP)$--generic (in the usual sense).  
\end{enumerate}
\item If $\bbP$ is proper over $(D,\cS)$--semi diamonds, $\mu\geq\lambda$, 
$Y\subseteq [\mu]^{\le\lambda}$, $Y\in\bV$, \underline{then}: 
\begin{enumerate}
\item[(a)] forcing with $\bbP$ does not collapse $\lambda^+$,
\item[(b)] forcing with $\bbP$ preserves the following two properties: 
\begin{enumerate}
\item[(i)]  $Y$ is a cofinal subset of $[\mu]^{\le \lambda}$ (under
inclusion), 
\item[(ii)] for every large enough $\chi$ and $x \in {\mathcal H}(\chi)$, 
there is $N\prec ({\mathcal H}(\chi),\in)$ such that $\|N\|=\lambda$, $N\cap
\lambda^+\in\lambda^+$, $N^{<\lambda}\subseteq N$, $N\cap\mu\in Y$ (i.e.,
the stationarity of $Y$ under the relevant filter).
\end{enumerate}
\end{enumerate}
\end{enumerate}
\end{proposition}

\begin{proof}
1)\quad Immediate (by the $\lambda$--completeness of $\bbP$).
\smallskip

\noindent 2)\quad Clause (a) should be clear (remember
\ref{pre1.2}(2)$(\alpha)$). For clause (b) note that $0 \in {\mathcal S}'$,
so in the game $\Game(r,N,h,\bbP,\bar{F},\bar{q})$ the condition $r_0$ is
chosen by the anti-generic player. So if the conclusion fails, then for some
$\bbP$--name $\name{\alpha}\in N$ for an ordinal we have $r\not\forces
\mbox{`` }\name{\alpha}\in N\mbox{ ''}$. Thus the anti-generic player can
choose $r_0$ so that $r_0\forces\mbox{`` }\name{\alpha}=\alpha_0\mbox{ ''}$
for some ordinal $\alpha_0\notin N$, what guarantees it to win the play. 
\smallskip

\noindent 3)\quad Straightforward from 2).
\end{proof}

Very often checking properness over $D$--semi diamonds (for particular
examples of forcing notions) we get somewhat stronger properties, which
motivate the following definition.
 
\begin{definition}
\label{gencon}
We say that a condition $r\in\bbP$ is {\em $N$--generic for $D$--semi
diamonds\/} if it is $(N,h,\bbP)$--generic for $\bar{q}$ over $\bar{F}$
whenever $h,\bar{q},\bar{F}$ are as in \ref{pre1.2}. Omitting $D$ we mean
``for every normal filter $D$ with diamonds''.  
\end{definition}

The following notion is not of main interest in this paper, but surely it is
interesting from the point of view of general theory.

\begin{definition}
\label{1.3} 
Let $0<\alpha<\lambda^+$.  
\begin{enumerate}
\item Let $\cS\in D^+$. We say that a forcing notion $\bbP$ is {\em
$\alpha$--proper over $(D,\cS)$--semi diamonds\/} whenever
\begin{enumerate}
\item[(a)] $\bbP$ is $\lambda$--complete, and 
\item[(b)] \underline{if} $\chi$ is large enough, $p\in\bbP$ and 
\begin{itemize} 
\item $\bar{N}=\langle N_\beta:\beta<\alpha\rangle$ is an increasing
sequence of elementary submodels of $({\mathcal H}(\chi),\in)$ such that
$\|N_\beta\|=\lambda$, $N^{< \lambda}_\beta \subseteq N_\beta$,
$\{\lambda,p,\bbP,\bar{N}\restriction \beta\}\in N_\beta$, and 
\item $\bar{F}^\beta=\langle F^\beta_\delta:\delta\in \cS\rangle$,
$F^\beta_\delta\in {}^\delta\delta$ (for $\beta<\alpha$),
\item $h_\beta:\lambda\longrightarrow N_\beta$, $\bbP\cap N_\beta\subseteq
\Rang(h_\beta)$ and $\langle h_\gamma,\bar{F}^\gamma:\gamma<\beta\rangle\in
N_\beta$, and  
\item $\bar{F}^\beta$ is a $(D,\cS)$--semi diamond sequence for $(N_\beta,
h_\beta,\bbP)$, and 
\item $\bar{q}^\beta=\langle q^\beta_\delta:\delta\in {\mathcal S}\rangle$
is an $(N_\beta,h_\beta,\bbP)$--candidate over $\bar{F}^\beta$, and $\langle
\bar{q}^\gamma:\gamma<\beta\rangle\in N_\beta$,  
\end{itemize}
\underline{then} there is $r\in\bbP$ above $p$ which is $(N_\beta,h_\beta,
\bbP)$--generic for $\bar{q}^\beta$ over $\bar{F}^\beta$ for each $\beta<
\alpha$.  
\end{enumerate}
\item We define ``{\em $\bbP$ is $\alpha$--proper over $D$--semi diamonds}''
(and $K^\alpha_D$) and ``{\em $\bbP$ is $\alpha$--proper over $\lambda$}''
in a way parallel to \ref{1.2}(2,3). 
\end{enumerate}
\end{definition}

\begin{remark}
Note that for $\alpha = 1$ (in Definition \ref{1.3}) we get the same notions
as in Definition \ref{1.2}.  
\end{remark}

\section{The preservation theorem}
In \ref{1.6}  below we prove a preservation theorem for our forcing
notions. It immediately gives the consistency of the suitable Forcing Axiom,
see \ref{axiom}. Also the proof actually specifies which semi-diamond
sequences $\bar{F}$ are used.  

First, recall that

\begin{proposition}
\label{first}
Suppose that $\langle \bbP_\alpha,\name{\bbQ}_\alpha:\alpha<\zeta^*\rangle$
is a $(< \lambda^+)$--support iteration such that for each $\alpha<\zeta^*$ 
\[\forces_{\bbP_\alpha}\mbox{`` $\name{\bbQ}_\alpha$ is
$\lambda$--complete. ''}\]
Then the forcing $\bbP_{\zeta^*}$ is $\lambda$--complete. 
\end{proposition}

Before we engage in the proof of the preservation theorem, let us prove some
facts of more general nature than the one of our main context. If, e.g., all
iterands are $\lambda$--lub--complete, then Proposition \ref{cl0} below is
obvious.  

\begin{tcontext}
Let $\bar{\bbQ}=\langle\bbP_\alpha,\name{\bbQ}_\alpha:\alpha<\zeta^*\rangle$
be a $(<\lambda^+)$--support iteration of $\lambda$--complete forcing
notions. We also suppose that $N$ is a model as in \ref{pre1.2},
$\bar{\bbQ},\ldots\in N$. 
\end{tcontext}

\begin{proposition}
\label{cl0}
Suppose that $\zeta\in(\zeta^*+1)\cap N$ is a limit ordinal of cofinality
$\cf(\zeta)<\lambda$ and $r\in\bbP_\zeta$ is such that 
\[(\forall\varepsilon\in\zeta\cap N)\big(r\restriction\varepsilon\mbox{ is
$(N,\bbP_\varepsilon)$--generic.}\big)\]
Assume that conditions $s_\beta\in N\cap\bbP_\zeta$ (for $\beta<\delta$,
$\delta<\lambda$) are such that  
\[(\forall\beta'<\beta<\delta)(s_{\beta'}\leq s_\beta\leq r).\]
\underline{Then} there are conditions $s\in N\cap\bbP_\zeta$ and $r^+\in
\bbP_\zeta$ such that $s\leq r^+$, $r\leq r^+$ and $(\forall\beta<\delta)
(s_\beta\leq s)$. 
\end{proposition}

\begin{proof}
Let $\langle i_\gamma:\gamma<\cf(\zeta)\rangle\subseteq N\cap\zeta$ be a
strictly increasing continuous sequence cofinal in $\zeta$. By induction on
$\gamma$ choose $r^-_\gamma,r_\gamma$ such that 
\begin{enumerate}
\item[$(\alpha)$] $r^-_\gamma\in \bbP_{i_\gamma}\cap N$ is above (in
$\bbP_{i_\gamma}$) of all $s_\beta\restriction i_\gamma$ for $\beta<\delta$, 
\item[$(\beta)$]  $r_\gamma\in\bbP_{i_\gamma}$, $r^-_\gamma
\leq_{\bbP_{i_\gamma}} r_\gamma$, and $r\restriction i_\gamma\leq r_\gamma$, 
\item[$(\gamma)$] if $\gamma<\varepsilon<\cf(\zeta)$ then $r^-_\gamma\leq
r^-_\varepsilon$ and $r_\gamma\leq r_\varepsilon$.
\end{enumerate}
(The choice is clearly possible as $r\restriction i_\gamma$ is $(N,
\bbP_{i_\gamma})$--generic.) 

Let $r^+\in\bbP_\zeta$ be an upper bound of $\langle r_\gamma:\gamma<\cf(
\zeta)\rangle$ (remember clause $(\gamma)$ above); then also $r\leq r^+$. 
Now we are going to define a condition $s\in\bbP_\zeta\cap N$. We let
$\Dom(s)=\bigcup\{\Dom(r^-_{\gamma+1})\cap [i_\gamma,i_{\gamma+1}):\gamma<
\cf(\zeta)\}$, and for $\xi\in\Dom(s)$, $i_\gamma\leq\xi<i_{\gamma+1}$, we
let $s(\xi)$ be a $\bbP_\xi$--name for the following object in
$\bV[G_{\bbP_\xi}]$ (for a generic filter $G_{\bbP_\xi}\subseteq\bbP_\xi$
over $\bV$):  
\begin{enumerate}
\item[(i)]   If $r_{\gamma+1}^-(\xi)[G_{\bbP_\xi}]$ is an upper bound of
$\{s_\beta(\xi)[G_{\bbP_\xi}]:\beta<\delta\}$ in $\name{\bbQ}_\xi[
G_{\bbP_\xi}]$,\\
\underline{then} $s(\xi)[G_{\bbP_\xi}]=r^-_{\gamma+1}(\xi)[G_{\bbP_\xi}]$.
\item[(ii)]  If not (i), but $\{s_\beta(\xi)[G_{\bbP_\xi}]:\beta<\delta\}$
has an upper bound in $\name{\bbQ}_\xi[G_{\bbP_\xi}]$,\\
\underline{then} $s(\xi)[G_{\bbP_\xi}]$ is the $<^*_\chi$--first such upper
bound. 
\item[(iii)] If neither (i) nor (ii), then $s(\xi)[G_{\bbP_\xi}]=s_0(\xi)
[G_{\bbP_\xi}]$. 
\end{enumerate}
It should be clear that $s\in\bbP_\zeta\cap N$. Now, 
\begin{itemize}
\item $s\leq r^+$.
\end{itemize}
Why? By induction on $\xi\in\zeta\cap N$ we show that $s\restriction\xi\leq
r^+ \restriction\xi$. Steps ``$\xi=0$'' and ``$\xi$ limit'' are clear, so
suppose that we have proved $s\restriction\xi\leq r^+\restriction\xi$,
$i_\gamma\leq\xi<i_{\gamma+1}$ (and we are interested in the restrictions to
$\xi+1$). Assume that $G_{\bbP_\xi}\subseteq\bbP_\xi$ is a generic filter
over $\bV$ such that $r^+\restriction\xi\in G_{\bbP_\xi}$. Since $s_\beta
\restriction i_{\gamma+1}\leq r^-_{\gamma+1}\leq r_{\gamma+1}\leq r^+$, we
also have $\{s_\beta\restriction\xi:\beta<\delta\}\subseteq G_{\bbP_\xi}$
and $r^-_{\gamma+1}\restriction\xi\in G_{\bbP_\xi}$. Hence $r^-_{\gamma+1}
(\xi)[G_{\bbP_\xi}]$ is an upper bound of $\{s_\beta(\xi)[G_{\bbP_\xi}]:
\beta<\delta\}$. Therefore, $s(\xi)[G_{\bbP_\xi}]=r^-_{\gamma+1}(\xi)[
G_{\bbP_\xi}]\leq r_{\gamma+1}(\xi)[G_{\bbP_\xi}]\leq r^+(\xi)[
G_{\bbP_\xi}]$ (see (i) above) and we are done.

The proof of the proposition will be finished once we show
\begin{itemize}
\item $(\forall\beta<\delta)(s_\beta\leq s)$.
\end{itemize}
Why does this hold? By induction on $\xi\in\zeta\cap N$ we show that
$s_\beta\restriction\xi\leq s\restriction\xi$ for all $\beta<\delta$. Steps
``$\xi=0$'' and ``$\xi$ limit'' are as usual clear, so suppose that we have
proved $s_\beta\restriction\xi\leq s\restriction\xi$ (for $\beta<\delta$),
$i_\gamma\leq\xi<i_{\gamma+1}$ (and we are interested in the restrictions to 
$\xi+1$). Assume that $G_{\bbP_\xi}\subseteq\bbP_\xi$ is a generic filter
over $\bV$ such that $s\restriction\xi\in G_{\bbP_\xi}$. Then also (by the
inductive hypothesis) $\{s_\beta\restriction\xi:\beta<\delta\}\subseteq
G_{\bbP_\xi}$ and therefore $\langle s_\beta(\xi)[G_{\bbP_\xi}]:\beta<\delta
\rangle$ is an increasing sequence of conditions from the
($\lambda$--complete) forcing $\name{\bbQ}_\xi[G_{\bbP_\xi}]$. Thus this
sequence has an upper bound, and $s(\xi)[G_{\bbP_\xi}]$ is such an upper
bound (see (i) and (ii) above), as required.
\end{proof}

In the proof of the preservation theorem we will (like in the proof of the
preservation of properness \cite[Ch.III, \S3.3]{Sh:f}) have to deal with
{\em names\/} for conditions in the iteration. This motivates the following
definition (which is in the spirit of \cite[Ch.X]{Sh:f}, so this is why
``RS'').

\begin{definition}
\label{RSconditions}
\begin{enumerate}
\item An {\em RS--condition in $\bbP_{\zeta^*}$} is a pair $(p,w)$ such that
$w\in [(\zeta^*+1)]^{<\lambda}$ is a closed set, $0,\zeta^*\in w$, $p$ is a
function with domain $\Dom(p)\subseteq \zeta^*$, and
\begin{enumerate}
\item[$(\otimes)_1$] for every two successive members $\varepsilon'<
\varepsilon''$ of the set $w$, $p\restriction [\varepsilon',\varepsilon'')$
is a $\bbP_{\varepsilon'}$--name of an element of $\bbP_{\varepsilon''}$
whose domain is included in the interval $[\varepsilon',\varepsilon'')$.
\end{enumerate}
The family of all RS--conditions in $\bbP_{\zeta^*}$ is denoted by
$\bbP_{\zeta^*}^{\rm RS}$.
\item If $(p,w)\in\bbP_{\zeta^*}^{\rm RS}$ and $G_{\bbP_{\zeta^*}}\subseteq
\bbP_{\zeta^*}$ is a generic filter over $\bV$, then we write $(p,w)\in'
G_{\bbP_{\zeta^*}}$ whenever
\begin{enumerate}
\item[$(\otimes)_2$] for every two successive members $\varepsilon'<
\varepsilon''$ of the set $w$, 
\[(p\restriction[\varepsilon',\varepsilon''))[G_{\bbP_{\zeta^*}}\cap
\bbP_{\varepsilon'}]\in G_{\bbP_{\zeta^*}}\cap\bbP_{\varepsilon''}.\]
\end{enumerate}
\item If $(p_1,w_1), (p_2,w_2)\in \bbP_{\zeta^*}^{\rm RS}$, then we write
$(p_1,w_1)\leq'(p_2,w_2)$ whenever 
\begin{enumerate}
\item[$(\otimes)_3$] for every generic $G_{\bbP_{\zeta^*}}\subseteq
\bbP_{\zeta^*}$ over $\bV$, if $(p_2,w_2)\in' G_{\bbP_{\zeta^*}}$ then
$(p_1,w_1)\in'G_{\bbP_{\zeta^*}}$ and for each two successive members
$\varepsilon'<\varepsilon''$ of the set $w_1\cup w_2$ we have
\[(p_1\restriction[\varepsilon',\varepsilon''))[G_{\bbP_{\zeta^*}}\cap
\bbP_{\varepsilon'}]\leq_{\bbP_{\varepsilon''}} (p_2\restriction[
\varepsilon',\varepsilon''))[G_{\bbP_{\zeta^*}}\cap\bbP_{\varepsilon'}].\] 
\end{enumerate}
\end{enumerate}
\end{definition}

\begin{remark}
If $(p,w)\in \bbP_{\zeta^*}^{\rm RS}$, $\varepsilon'\leq\xi<\varepsilon''$,
$\varepsilon',\varepsilon''$ are successive members of $w$, then $p(\xi)$ is
a $\bbP_{\varepsilon'}$--name for a $\bbP_\xi$--name of a member of
$\name{\bbQ}_\xi$. One may look at this name as a $\bbP_\xi$--name. However
note that if we apply this approach to each $\xi$, we may not end up with a
condition in $\bbP_{\zeta^*}$ because of the support!
\end{remark}

\begin{proposition}
\label{RSbounds}
\begin{enumerate}
\item For each $(p,w)\in\bbP_{\zeta^*}^{\rm RS}$ there is $q\in
\bbP_{\zeta^*}$ such that $(p,w)\leq'(q,\{0, \zeta^*\})$. 
\item If $(p,w)\in\bbP_{\zeta^*}$ and $q\in\bbP_{\zeta^*}$, then there is
$q^*\in\bbP_{\zeta^*}$ stronger than $q$ and such that for each successive
members $\varepsilon'<\varepsilon''$ of $w$ the condition $q^*\restriction
\varepsilon'$ decides $p\restriction [\varepsilon',\varepsilon'')$ (i.e.,
$q\restriction\varepsilon'\forces$`` $p\restriction [\varepsilon',
\varepsilon'')=p_{\varepsilon',\varepsilon''}$ '' for some $p_{\varepsilon',
\varepsilon''}\in\bbP_{\zeta^*}$).
\item Let $(p_i,w_i)\in \bbP_{\zeta^*}^{\rm RS}\cap N$ (for $i<\delta<
\lambda$), and $s\in\bbP_{\zeta^*}\cap N$, $r\in\bbP_{\zeta^*}$ be such that 
\[s\leq r\quad\mbox{ and }\quad(\forall j<i<\delta)((p_j,w_j)\leq' (p_i,w_i)
\leq' (r,\{0,\zeta^*\}).\]
Assume that either $r$ is $(N,\bbP_{\zeta^*})$--generic, or $\zeta^*$ is a
limit ordinal of cofinality $\cf(\zeta^*)<\lambda$ and for every $\zeta<
\zeta^*$ the condition $r\restriction\zeta$ is $(N,\bbP_\zeta)$--generic.

\noindent \underline{Then} there are conditions $s'\in N\cap\bbP_{\zeta^*}$
and $r'\in\bbP_{\zeta^*}$ such that $s\leq s'\leq r'$, $r\leq r'$ and
$(\forall i<\delta)((p_i,w_i)\leq' (s',\{0,\zeta^*\}))$.  
\end{enumerate}
\end{proposition}

\begin{proof}
(1), (2)\qquad Straightforward (use the $\lambda$--completeness of
$\bbP_{\zeta^*}$).  
\medskip

\noindent (3)\quad If $r$ is $(N,\bbP_{\zeta^*})$--generic, then our
assertion is clear (remember clause (2)). So suppose that we are in the
second case (so $\aleph_0\leq\cf(\zeta^*)<\lambda$). Let $\langle i_\gamma:
\gamma<\cf(\zeta^*)\rangle\subseteq N\cap\zeta$ be a strictly increasing
continuous sequence cofinal in $\zeta^*$. For $\gamma<\cf(\zeta^*)$ and
$i<\delta$ let $p_i^\gamma=p_i\restriction i_\gamma$, $w^\gamma_i=(w_i\cap
i_\gamma)\cup\{i_\gamma\}$ (clearly $(p^\gamma_i,w^\gamma_i)\in
\bbP_{i_\gamma}^{\rm RS}$). Since $r\restriction i_\gamma$ is $(N,
\bbP_{i_\gamma})$--generic, we may inductively pick conditions $s_\gamma,
r_\gamma$ (for $\gamma<\cf(\zeta^*)$) such that
\begin{itemize}
\item $s\restriction i_\gamma\leq s_\gamma\in\bbP_{i_\gamma}\cap N$, $r\leq
r_\gamma\in\bbP_{\zeta^*}$,
\item $(\forall i<\delta)((p^\gamma_i,w^\gamma_i)\leq' (s_\gamma,
\{0,i_\gamma\}))$, $s_\gamma\leq r_\gamma\restriction i_\gamma$,  
\item if $\beta<\gamma<\cf(\zeta^*)$ then $s\restriction i_\beta\leq s_\beta
\leq s_\gamma$ and $r_\beta\leq r_\gamma$.
\end{itemize}
Let $r^*\in \bbP_{\zeta^*}$ be stronger than all $r_\gamma$'s. Now apply
\ref{cl0}. 
\end{proof}
\medskip

Now we may state and prove our main result.

\begin{theorem}
\label{1.6} 
Let $D,{\mathcal S},{\mathcal S}'$ be as in \ref{1.1}. Assume that
$\bar{\bbQ}=\langle\bbP_\alpha,\name{\bbQ}_\alpha:\alpha<\zeta^*\rangle$ is
a $(<\lambda^+)$--support iteration such that for each $\alpha<\zeta^*$ 
\[\forces_{\bbP_\alpha}\mbox{`` $\name{\bbQ}_\alpha$ is proper for
$D$--semi diamonds ''.}\] 
\underline{Then} $\bbP_{\zeta^*}=\lim(\bar{\bbQ})$ is proper for
$D$--semi diamonds.  
\end{theorem}

\begin{proof}  
By \ref{first}, the forcing notion $\bbP_{\zeta^*}$ is $\lambda$--complete,
so we have to concentrate on showing clause \ref{1.2}(1)(b) for it.  

So suppose that $\chi$ is large enough, $p\in\bbP_{\zeta^*}$ and $N\prec
({\mathcal H} (\chi),\in,<^*_\chi)$, $\|N\|=\lambda$, $N^{<\lambda}\subseteq 
N$ and $\{\lambda,p,\bar{\bbQ},\bbP_{\zeta^*},D,\cS,\ldots\}\in N$, and
$h:\lambda\longrightarrow N$  satisfies $\bbP_{\zeta^*}\cap N\subseteq
\Rang(h)$. Furthermore, suppose that $\bar{F}=\langle F_\delta:\delta\in
\cS\rangle$ is a $(D,\cS)$--semi diamond and $\bar{q}=\langle q_\delta:
\delta\in\cS\rangle$ is an $(N,h,\bbP_{\zeta^*},\bar{F})$--candidate. 
We may assume that for each $\delta\in\cS$ 
\begin{enumerate}
\item[$(\odot)$] \underline{if} $\langle h\circ F_\delta(\alpha):\alpha<
\delta\rangle$ is not a $\leq_{\bbP_{\zeta^*}}$--increasing sequence of
members of $\bbP_{\zeta^*}\cap N$, 

\noindent \underline{then} $h\circ F_\delta(\alpha)=*$ for all $\alpha<
\delta$. 
\end{enumerate}
[Just suitably modify $F_\delta$ whenever the assumption of $(\odot)$ holds
-- note that the modification does not change the notion of a candidate, the
game from \ref{pre1.2}(3), etc.]  
\medskip

Before we define a generic condition $r\in\bbP_{\zeta^*}$ for $\bar{q}$ over
$\bar{F}$, let us introduce notation used later and give two important
facts.  
\medskip

Let $i\in N\cap(\zeta^*+1)$ and let $G_{\bbP_i}\subseteq\bbP_i$ be generic 
over $\bV$. We define:
\begin{itemize}
\item $h^{\langle i\rangle}[G_{\bbP_i}]:\lambda\longrightarrow
N[G_{\bbP_i}]$ is such that

\noindent if $h(\gamma)$ is a function, $i\in\Dom(h(\gamma))$ and $(h(
\gamma))(i)$ is a $\bbP_i$--name, then $(h^{\langle i\rangle}[G_{\bbP_i}])
(\gamma)=(h(\gamma))(i)[G_{\bbP_i}]$, otherwise it is $*$; 
\item $h^{[i]}:\lambda\longrightarrow N$ is defined by 

\noindent $h^{[i]}(\gamma)= (h(\gamma))\restriction i$ provided $h(\gamma)$
is a function, and $*$ otherwise;
\item ${\mathcal S}^{\langle i\rangle}[G_{\bbP_i}]=\{\delta\in {\mathcal
S}:\mbox{ if $\delta$ is limit, then }q_\delta \restriction i\in
G_{\bbP_i}\}$; 
\item $\bar{q}^{\langle i\rangle}[G_{\bbP_i}]$ is $\langle q_\delta(i)[
G_{\bbP_i}]:\delta\in{\mathcal S}^{\langle i\rangle}[G_{\bbP_i}]\rangle$; 
\item $\bar{q}^{[i]}=\langle q_\delta\restriction i:\delta\in{\mathcal S}
\rangle$; 
\item $\bar{F}^{\langle i\rangle}[G_{\bbP_i}]$ is $\langle F_\delta:\delta
\in {\mathcal S}^{\langle i\rangle}[G_{\bbP_i}]\rangle$.
\end{itemize}
Observe that $h^{[i]}:\lambda\longrightarrow N$ is such that $\bbP_i\cap N
\subseteq\Rang(h^{[i]})$ and $h^{\langle i\rangle}[G_{\bbP_i}]$ is such that
$N[G_{\bbP_i}]\cap \name{\bbQ}_i[G_{\bbP_i}]\subseteq\Rang(h^{\langle
i\rangle}[G_{\bbP_i}])$. 

Plainly, by $(\odot)$,  
\begin{claim}
\label{cl1}
Assume $i\in N\cap (\zeta^*+1)$. \underline{Then} $\bar{F}$ is a
$(D,\cS)$--semi diamond sequence for $(N,h^{[i]},\bbP_i)$ and
$\bar{q}^{[i]}$ is an $(N,h^{[i]},\bbP_i,\bar{F})$--candidate.  
\end{claim}

\begin{claim}
\label{cl2}
Assume that $i\in N\cap (\zeta^*+1)$ and $r\in\bbP_i$ is $(N,h^{[i]},
\bbP_i)$--generic for $q^{[i]}$ over $\bar{F}$. Let $G_{\bbP_i}\subseteq
\bbP_i$ be a generic filter over $\bV$, $r\in G_{\bbP_i}$. \underline{Then}
in $\bV[G_{\bbP_i}]$: 
\begin{enumerate}
\item ${\mathcal S}^{\langle i\rangle}[G_{\bbP_i}]\in D^+$,
\item $\bar{F}^{\langle i\rangle}[G_{\bbP_i}]$ is a
$(D,\cS^{\langle i\rangle}[G_{\bbP_i}])$--semi diamond for 

$(N[G_{\bbP_i}],h^{\langle i\rangle}[G_{\bbP_i}],\name{\bbQ}_i[
G_{\bbP_i}])$, and 
\item $\bar{q}^{\langle i\rangle}[G_{\bbP_i}]$ is an
$(N[G_{\bbP_i}],h^{\langle i\rangle}[G_{\bbP_i}],\name{\bbQ}_i[G_{\bbP_i}],
\bar{F}^{\langle i\rangle}[G_{\bbP_i}])$--candidate.  
\end{enumerate}
\end{claim}

\begin{proof}[Proof of the Claim]
(1)\quad Will follow from (2).
\medskip

\noindent (2)\quad Assume that this fails. Then we can find a condition
$r^*\in\bbP_i$, a $\bbP_i$--name $\name{\bar{q}}'=\langle\name{q}_\alpha':
\alpha<\lambda\rangle\subseteq N$ for an increasing sequence of conditions
from $\name{\bbQ}_i$, and $\bbP_i$--names $\name{A}_\xi$ for members of
$D\cap\bV$ such that $r\leq_{\bbP_i} r^*\in G_{\bbP_i}$ and  
\[r^*\forces_{\bbP_i}\mbox{`` }(\forall\delta\in\cS^{\langle i\rangle}\cap 
\mathop{\triangle}\limits_{\xi<\lambda}\name{A}_\xi)(\langle h^{\langle
i\rangle}\circ F_\delta(\alpha):\alpha<\delta\rangle\neq\name{\bar{q}}'
\restriction \delta)\mbox{ ''.}\]  
Consider a play $\langle r_j^-,r_j,C_j:i<\lambda\rangle\subseteq\bbP_i$ of
the game $\Game(r,N,h^{[i]},\bbP_i,\bar{F},\bar{q}^{[i]})$ in which the
generic player uses its winning strategy and the anti-generic player plays
as follows. In addition to keeping the rules of the game, it makes sure that 
at stage $j\in\cS'$: 
\begin{itemize}
\item $r_j\geq r^*$ (so $r_0\geq r^*$; remember the anti-generic player plays
at 0), 
\item $r_j$ decides the values of all $\name{A}_\xi$ for $\xi<j$. 
\end{itemize}
Let $A_\xi\in D\cap\bV$ be such that $r_j\forces\mbox{`` }\name{A}_\xi=
A_\xi\mbox{ ''}$ for sufficiently large $j\in {\mathcal S}'$. 

Note that the sequence $\langle {r_j^-}^\frown\langle\name{q}_j'\rangle:j<
\delta\rangle$ is $\leq_{\bbP_{i+1}}$--increasing. So, as $D$ is normal and
$A_\xi,C_j\in D$ and $\bar{F}$ is a semi diamond for $(N,h^{[i+1]},
\bbP_{i+1})$ (by \ref{cl1}), we may find a limit ordinal $\delta\in\cS\cap
\mathop{\triangle}\limits_{\xi<\lambda}A_\xi\cap\mathop{\triangle}\limits_{
j<\lambda}C_j$ such that $\langle h^{[i+1]}\circ F_\delta(j):j<\delta\rangle
=\langle {r_j^-}^\frown\langle\name{q}_j'\rangle:j<\delta\rangle$. Then
also $\langle h^{[i]}\circ F_\delta(j):j<\delta\rangle=\langle r_j^-:j<
\delta\rangle$, and since the play is won by the generic player, we conclude
$q_\delta\restriction i\leq r_\delta$. But then taking sufficiently large
$j\in\cS'$ we have 
\[r_j\forces\mbox{`` }\delta\in {\mathcal S}^{\langle i\rangle}\cap
\mathop{\triangle}_{\xi<\lambda}\name{A}_\xi\quad \&\quad \langle h^{\langle
i\rangle}\circ F_\delta(\alpha):\alpha<\delta\rangle=\name{\bar{q}}'
\restriction\delta\mbox { '',}\] 
a contradiction.
\medskip

\noindent (3)\quad Should be clear. 
\end{proof}
\medskip

Fix a bijection $\Upsilon:\zeta^*\cap N\longrightarrow\gamma^*\leq\lambda$. 
Also let $\langle(\name{\tau}_i,\zeta_i):i<\lambda\rangle$ list all pairs
$(\name{\tau},\zeta)\in N$ such that $\zeta\leq \zeta^*$, $\cf(\zeta)\geq
\lambda$ and $\name{\tau}$ is a $\bbP_\zeta$--name for an ordinal.  

Next, by induction, we choose a sequence $\langle (p_i,w_i):i<\lambda
\rangle\subseteq\bbP_{\zeta^*}^{\rm RS}\cap N$ such that 
\begin{enumerate}
\item[(i)]   $(p,\{0,\zeta^*\})\leq' (p_i,w_i)\leq' (p_j,w_j)$ for $i<j<
\lambda$,
\item[(ii)]  if $i<j<\lambda$ and $\Upsilon(\varepsilon)\leq i$, then
$\varepsilon\in\Dom(p_i)$ and $p_i(\varepsilon)=p_j(\varepsilon)$,
\item[(iii)] if $i<\lambda$ is a limit ordinal, then $w_i$ is the closure of
$\bigcup\limits_{j<i}w_j$, and\\
if, additionally, $\varepsilon\in\Dom(q_i)$ is such that $\Upsilon(
\varepsilon)\geq i$ (and $i\in\cS$, of course), then $\varepsilon\in
\Dom(p_i)$ and $p_i(\varepsilon)$ is such that 
\begin{enumerate}
\item[$(\otimes)$] for every generic $G_{\bbP_{\zeta^*}}\subseteq
\bbP_{\zeta^*}$ over $\bV$ such that $(p_i,w_i)\in' G_{\bbP_{\zeta^*}}$,  
and two successive members $\varepsilon',\varepsilon''$ of the set $w_i$
such that $\varepsilon'\leq \varepsilon<\varepsilon''$ we have:

\underline {if} $\{p_j(\varepsilon)[G_{\bbP_{\zeta^*}}\cap
\bbP_{\varepsilon'}][G_{\bbP_{\zeta^*}}\cap \bbP_\varepsilon]:j<i\}\cup
\{q_i(\varepsilon)[G_{\bbP_{\zeta^*}}\cap \bbP_\varepsilon]\}$ has an upper 
bound in $\name{\bbQ}_\varepsilon[G_{\bbP_{\zeta^*}}\cap\bbP_\varepsilon
]$,\\  
\underline{then} $p_i(\varepsilon)[G_{\bbP_{\zeta^*}}\cap\bbP_{\varepsilon'}
][G_{\bbP_{\zeta^*}}\cap \bbP_\varepsilon]$ is such an upper bound,
\end{enumerate}
\item[(iv)] for each $i<\lambda$, for some $\xi\in N\cap\zeta_i$ and a
$\bbP_\xi$--name $\name{\tau}\in N$ we have:\\
$\sup\big(\{\varepsilon<\zeta_i:\Upsilon(\varepsilon)\leq i\}\cup
(w_i\cap\zeta_i)\big)<\xi$, $w_{i+1}=w_i\cup\{\xi\}$, $p_{i+1}\restriction
\xi=p_i\restriction \xi$ and\\ 
\underline{if} $G_{\bbP_{\zeta^*}}\subseteq\bbP_{\zeta^*}$ is generic over
$\bV$ and $(p_{i+1},w_{i+1})\in' G_{\bbP_{\zeta^*}}$,\\
\underline{then} $\name{\tau}_i[G_{\bbP_{\zeta^*}}\cap\bbP_{\zeta_i}]=
\name{\tau}[G_{\bbP_{\zeta^*}}\cap\bbP_\xi]$
\end{enumerate}
(It should be clear that there are no problems in the induction and it is
possible to pick $(p_i,w_i)$ as above.) From now on we will treat each
$p_i(\xi)$ as a $\bbP_\xi$--name for a member of $\name{\bbQ}_\xi$.

Now we are going to define an $(N,h,\bbP_{\zeta^*})$--generic condition 
$r\in\bbP$ for $\bar{q}$ over $\bar{F}$ in the most natural way. Its domain
is $\Dom(r)=\zeta^*\cap N$ and for each $i\in\zeta^*\cap N$
\[r\restriction i\forces\mbox{`` }r(i)\geq p_{\Upsilon(i)}(i)\mbox{ is $(N[ 
\name{G}_{\bbP_i}], h^{\langle i\rangle},\name{\bbQ}_i)$--generic for
$\bar{q}^{\langle i\rangle}$ over $\bar{F}^{\langle i\rangle}$ ''.}\]

\begin{mclaim}
\label{cl3}
For every $\zeta\in (\zeta^*+1)\cap N$, the generic player has a winning
strategy in the game $\Game(r\restriction\zeta,N,h^{[\zeta]},\bbP_\zeta,
\bar{F},\bar{q}^{[\zeta]})$. 
\end{mclaim}

\begin{proof}[Proof of the Claim]
We will prove the claim by induction on $\zeta\in (\zeta^*+1)\cap N$. For
$\zeta\in\zeta^*\cap N$ this implies that $r(\zeta)$ is well-defined
(remember \ref{cl2}). Of course for $\zeta=\zeta^*$ we finish the proof of
the theorem.     
\medskip

Suppose that $\zeta\in (\zeta^*+1)\cap N$ and we know that $r\restriction
\zeta'$ is $(N,h^{[\zeta']},\bbP_{\zeta'})$--generic for
$\bar{q}^{[\zeta']}$ over $\bar{F}$ for all $\zeta'\in N\cap\zeta$. We are
going to describe a winning strategy for the generic player in the game
$\Game(r\restriction\zeta,N,h^{[\zeta]},\bbP_\zeta,\bar{F},
\bar{q}^{[\zeta]})$. The inductive hypothesis is not used in the full
strength in the definition of the strategy, but we need it in several
places, e.g., to know that $r$ is well defined as well as that we have the
$\st_i$'s below. Also note that it implies that $(p_i^\zeta,w^\zeta_i)\leq'
(r\restriction\zeta,\{0,\zeta\})$ for all $i<\lambda$, where $p^\zeta_i=p_i
\restriction\zeta$ and $w^\zeta_i=(w_i\cap\zeta)\cup\{\zeta\}$. Moreover,
during the play, both players will always have legal moves. Why? By the
inductive hypothesis we know that $r\restriction\zeta'$ is $(N,
\bbP_{\zeta'})$--generic for all $\zeta'\in\zeta\cap N$. Therefore, if
$\zeta$ is a successor or a limit ordinal of cofinality $\geq\lambda$, then
we immediately get that $r\restriction\zeta$ is $(N,\bbP_\zeta)$--generic
(remember clause (iv) of the choice of the $p_i$'s!), and thus \ref{obs}(4)
applies. If $\zeta$ is a limit ordinal of cofinality $\cf(\zeta)<\lambda$,
then we may use \ref{cl0}. 
\medskip

Let $\st_i$ be a $\bbP_i$--name for the winning strategy of the generic
player in $\Game(r(i),N[\name{G}_{\bbP_i}],h^{\langle i\rangle},
\name{\bbQ}_i,\bar{F}^{\langle i\rangle},\bar{q}^{\langle i\rangle})$, and
let  
\[E_0\stackrel{\rm def}{=}\{\delta<\lambda:\delta\mbox{ is a limit of points
from $\cS'$ }\}.\]  
Plainly, $E_0$ is a club of $\lambda$.

Let the generic player play as follows. Aside it will construct sequences
$\langle \name{r}_{j'}^\ominus(\varepsilon),\name{r}_{j'}^\oplus(
\varepsilon):j'<\lambda,\ \varepsilon\in\zeta\cap N\rangle$ and $\langle
\name{C}_{j'}^\xi(\varepsilon):j',\xi<\lambda,\ \varepsilon \in\zeta\cap
N\rangle$ so that  
\begin{itemize}
\item $\name{r}^\ominus_{j'}(\varepsilon)$ is a $\bbP_\varepsilon$--name for
a member of $\name{\bbQ}_\varepsilon\cap N[\name{G}_{\bbP_\varepsilon}]$,
$\name{r}_{j'}^\oplus(\varepsilon)$ is a $\bbP_\varepsilon$--name for a 
member of $\name{\bbQ}_\varepsilon$, $\name{C}_{j'}^\xi(\varepsilon)$ is a 
$\bbP_\varepsilon$--name for a member of $D\cap\bV$, and
\item if $j\in\cS$, $j'\leq j$, and $\Upsilon(\varepsilon)\leq j$, then
after the $j^{\rm th}$ move (which is a move of the generic player) the
terms $\langle\name{C}_{j'}^\xi(\varepsilon):\xi<\lambda\rangle,
\name{r}_{j'}^\ominus(\varepsilon)$, and $\name{r}_{j'}^\oplus(\varepsilon)$
are defined.    
\end{itemize}
So suppose that $j^*\in\cS$ and $\langle r^-_j,r_j,C_j:j<j^*\rangle$ is the
result of the play so far. To clearly describe the answer of the generic
player we will consider two (only slightly different) cases in the order in
which they appear in the game. (Remember $(r_0,C_0)$ is chosen by the
anti-generic player and that all successor moves are done by the generic
player.) 
\medskip

\noindent{\sc Case 1:}\quad $j_0<j'<\min(\cS'\setminus (j_0+1))=j_1$, 
$j_0\in\cS'$.\\
First the generic player picks conditions $s^-,s\in\bbP_\zeta$, $s^-\in N$
such that $r_{j_0}^-\leq s^-\leq s$, $r_{j_0}\leq s$ and for each
$\xi\in\zeta\cap N$ we have
\[s^-\restriction\xi\forces\mbox{`` }(\forall i<j_0)(p_i(\xi)\leq s^-(\xi)) 
\mbox{ ''.}\]
[Why possible? By \ref{RSbounds}(3).]\\ 
Now the generic player looks at $\varepsilon_\gamma<\zeta$ such that
$\Upsilon(\varepsilon_\gamma)=\gamma<j_1$. It picks $\bbP_{ 
\varepsilon_\gamma}$--names $\name{r}_{j'}^\ominus(\varepsilon_\gamma),
\name{r}_{j'}^\oplus(\varepsilon_\gamma),\name{C}_{j'}^\xi(
\varepsilon_\gamma)$ so that $s\restriction \varepsilon_\gamma$ forces that  
\[\langle \name{r}_{j'}^\ominus(\varepsilon_\gamma),\name{r}_{j'}^\oplus(
\varepsilon_\gamma),\mathop{\triangle}\limits_{\xi<\lambda}\name{C}_{j'}^\xi
(\varepsilon_\gamma):j'<j_1\rangle\]
is a play according to $\st_{\varepsilon_\gamma}$ in which the moves of the
anti-generic player are determined as follows. First, it keeps the
convention that if $j'\in\cS\setminus\cS^{\langle\varepsilon_\gamma
\rangle}$, then $(\name{r}_{j'}^\ominus(\varepsilon_\gamma), 
\name{r}_{j'}^\oplus(\varepsilon_\gamma),\mathop{\triangle}\limits_
{ \xi<\lambda}\name{C}_{j'}^\xi(\varepsilon_\gamma))$ is (a name for) the
$<^*_\chi$--first legal answer to the play so far.  Now, if $\gamma<j_0$, 
then we have already the play up to $j_0$ (it easily follows from the
inductive construction that $s\restriction\varepsilon_\gamma$ indeed forces
that it is a ``legal'' play). The $j_0^{\rm th}$ move of the anti-generic
player is stipulated as $\name{r}_{j_0}^\ominus(\varepsilon_\gamma)=
s^-(\varepsilon_\gamma)$, $\name{r}_{j_0}^\oplus(\varepsilon_\gamma)=
s(\varepsilon_\gamma)$, $\name{C}_{j_0}^\xi(\varepsilon_\gamma)=
\bigcap\limits_{j\leq j_0}C_j$, and next we continue up to $j_1$ keeping our
convention. If $j_0\leq\gamma<j_1$, then the generic player lets
$\name{r}_0^\ominus(\varepsilon_\gamma)=s^-(\varepsilon_\gamma),
\name{r}_0^\oplus(\varepsilon_\gamma)=s(\varepsilon_\gamma)$, 
$\name{C}_0^\xi(\varepsilon_\gamma)=\bigcap\limits_{j\leq j_0}C_j$ and then
it ``plays'' the game according to $\st_{\varepsilon_\gamma}$ up to $j_1$
keeping our convention for all $j'\notin\cS^{\langle\varepsilon_\gamma
\rangle}$.

Next, the generic player picks a condition $r^*\in\bbP_{\zeta}$ and
$\bbP_{\varepsilon_\gamma}$--names $\name{\tau}_{j'}(\varepsilon_\gamma)\in
N$ (for $\gamma<j_1$, $\varepsilon_\gamma<\zeta$, $j_0<j'<j_1$) such that
\begin{itemize}
\item $r^*\geq s$, and for every $\gamma,j'<j_1$, and
\[r^*\restriction\varepsilon_\gamma\forces_{\bbP_{\varepsilon_\gamma}}
\mbox{`` }\name{r}_{j'}^\oplus(\varepsilon_\gamma)\leq r^*(
\varepsilon_\gamma)\ \ \&\ \ \name{r}_{j'}^\ominus(\varepsilon_\gamma)=
\name{\tau}_{j'}(\varepsilon_\gamma)\mbox{ ''},\] 
\item for every $j,\xi<j_1$ and $\gamma<j_1$ with $\varepsilon_\gamma<
\zeta$, the condition $r^*\restriction\varepsilon_\gamma$ decides the value
of $\name{C}_{j'}^\xi(\varepsilon_\gamma)$, and 
\[r^*\restriction\varepsilon_\gamma\forces\mbox{`` }\name{C}_{j'}^\xi(
\varepsilon_\gamma)\setminus (\xi+1)=C_{j'}^\xi(\varepsilon_\gamma)
\mbox{ ''},\]
where $C_{j'}^\xi(\varepsilon_\gamma)\in D\cap\bV$.
\end{itemize}
Then it lets $r^-_{j'}\in N\cap\bbP_\zeta$ (for $j'\in (j_0,j_1)$) be
conditions such that 
\[\Dom(r^-_{j'})=\Dom(s^-)\cup\{\varepsilon_\gamma:\gamma<j_1\ \&\
\varepsilon_\gamma<\zeta\},\]
and for $\xi\in \Dom(r^-_{j'})$
\[\begin{array}{ll}
r^-_{j'}\restriction\xi\forces&\mbox{`` if $\Upsilon(\xi)<j_1$ and
$\langle\name{\tau}_j(\xi):j_0<j<j_1\rangle$ is an increasing}\\
&\mbox{sequence of conditions stronger than }s^-(\xi),\\
&\mbox{then }r_{j'}^-(\xi)=\name{\tau}_{j'} (\xi),\quad\mbox{ otherwise }
r_{j'}^-(\xi)=s^-(\xi)\mbox{ ''.} 
  \end{array}\]
Finally, for $j'\in (j_0,j_1)$ it plays $r^-_{j'},r^*,\bigcap\{C_{j'}^\xi( 
\varepsilon_\gamma):j',\gamma,\xi<j_1,\ \varepsilon_\gamma<\zeta\}$. 
\medskip

\noindent{\sc Case 2:}\quad $\sup\{i\in\cS':i<j'\}=j_0\leq j'<\min(\cS'
\setminus j_0)=j_1$, $j_0\in\cS$.\\
The generic player proceeds as above, the difference is that now $j_0$ 
``belongs to'' the generic player, and that it is a limit of moves of the
anti-generic player. Again, we look at $\varepsilon_\gamma<\zeta$ such that
$\Upsilon(\varepsilon_\gamma)=\gamma<j_1$.  

If $\gamma<j_0$, then every condition in $\bbP_{\varepsilon_\gamma}$
stronger than all $r_j\restriction\varepsilon_\gamma$ (for $j<j_0$) forces
that 
\[\langle \name{r}_{j'}^\ominus(\varepsilon_\gamma),\name{r}_{j'}^\oplus(
\varepsilon_\gamma),\mathop{\triangle}\limits_{\xi<\lambda}\name{C}_{j'}^\xi
(\varepsilon_\gamma):j'<j_0\rangle\]
is a legal play in which the generic player uses $\st_{
\varepsilon_\gamma}$. The generic player determines $\name{r}_{j'}^\ominus(
\varepsilon_\gamma),\name{r}_{j'}^\oplus(\varepsilon_\gamma)$, and
$\name{C}_{j'}^\xi(\varepsilon_\gamma)$ for $j'\in [j_0,j_1)$ ``playing the
game'' as earlier (with the same convention that if $j'\in\cS\setminus
\cS^{\langle\varepsilon_\gamma\rangle}$, then the $j'$--th move of the
anti-generic player is stipulated as the $<^*_\chi$--first legal move).

If $j_0\leq\gamma<j_1$, then (any condition stronger than all
$r_j\restriction\varepsilon_\gamma$ for $j<j_0$ forces that) $\langle
r^-_j(\varepsilon):j<j_0\rangle$, $\langle r_j(\varepsilon):j<j_0\rangle$
are increasing, and $r^-_j(\varepsilon_\gamma)\leq r_j(\varepsilon_\gamma)$
and $r(\varepsilon_\gamma)$ is $(N[\name{G}_{\bbP_{\varepsilon_\gamma}}],
\name{\bbQ}_{\varepsilon_\gamma})$--generic. So, by \ref{obs}(4), the
generic player may let $(r^\ominus_0(\varepsilon_\gamma),r^\oplus_0
( \varepsilon_\gamma))$ be the $<^*_\chi$--first such that for all $j<j_0$
we have $r^-_j(\varepsilon)\leq r^\ominus_0(\varepsilon_\gamma)\in N
[ \name{G}_{\bbP_{\varepsilon_\gamma}}]$, $r_j(\varepsilon_\gamma)\leq
r^\oplus_0(\varepsilon_\gamma)$. It also lets $C_0^\xi(\varepsilon_\gamma)
=\bigcap\limits_{j<j_0} C_j$. Then the generic player chooses
$\name{r}_{j'}^\ominus(\varepsilon_\gamma),\name{r}_{j'}^\oplus(
\varepsilon_\gamma)$, and $\name{C}_{j'}^\xi(\varepsilon_\gamma)$ for
$0<j'<j_1$ ``playing the game'' with the strategy $\st_{\varepsilon_\gamma}$
(and keeping the old convention for $j'\notin\cS^{\langle\varepsilon_\gamma
\rangle}$).  

Next the generic player picks a condition $r^*\in\bbP_{\zeta}$ (stronger
than all $r_j$ for $j<j_0$), $\bbP_{\varepsilon_\gamma}$--names
$\name{\tau}_{j'}(\varepsilon_\gamma)\in N$ and sets $C_{j'}^\xi(
\varepsilon_\gamma)\in D\cap\bV$ (for $j',\gamma,\xi<j_1$) as in the
previous case. Then it chooses conditions $s^-\in N\cap\bbP_\zeta$ and $r^+
\in\bbP_\zeta$ such that $r^*\leq r^+$ and $(\forall j<j_0)(r^-_j\leq s^-
\leq r^+)$. [Why possible? If $\zeta$ is limit of cofinality $\cf(\xi)<
\lambda$, use \ref{cl0}; otherwise we know that $r$ is
$(N,\bbP_\zeta)$--generic.] Next it defines conditions $r_{j'}^-\in N\cap
\bbP_\zeta$ (for $j_0\leq j'<j_1$) so that 
\[\Dom(r^-_{j'})=\Dom(s^-)\cup\{\varepsilon_\gamma:\gamma<j_1\ \&\
\varepsilon_\gamma<\zeta\},\]
and for $\xi\in \Dom(r^-_{j'})$
\[\begin{array}{ll}
r^-_{j'}\restriction\xi\forces&\mbox{`` if $\Upsilon(\xi)<j_1$ and
$\langle\name{\tau}_j(\xi):j_0\leq j<j_1\rangle$ is an increasing}\\
&\mbox{sequence of conditions above all $r^-_j(\xi)$ for $j<j_0$},\\
&\mbox{then }r_{j'}^-(\xi)=\name{\tau}_{j'} (\xi),\quad\mbox{ otherwise }
r_{j'}^-(\xi)=s^-(\xi)\mbox{ ''.} 
  \end{array}\]
Finally, for $j_0\leq j'<j_1$ it plays $r^-_{j'},r^+,\bigcap\{C_{j'}^\xi(
\varepsilon_\gamma):j',\gamma,\xi<j_1,\ \varepsilon_\gamma<\zeta\}$.  
\medskip

Why does the strategy described above work? Suppose that $\langle r_j,C_j:j<
\lambda\rangle$ is a play of the game  $\Game(r\restriction\zeta,N,
h^{[\zeta]},\bbP_{\zeta},\bar{F},\bar{q}^{[\zeta]})$ in which the generic
player used this strategy and let $\langle \name{r}_{j'}'(\varepsilon):j'<
\lambda,\ \varepsilon\in\zeta\cap N\rangle$ and $\langle\name{C}_{j'}^\xi(
\varepsilon):j',\xi<\lambda,\ \varepsilon \in\zeta\cap N\rangle$ be the
sequences it constructed aside. (As we said earlier, the game surely lasted
$\lambda$ steps and thus the sequences described above have length
$\lambda$.) 

Let us argue that condition \ref{pre1.2}(3)$(\circledast)$ holds.
\medskip

Assume that a limit ordinal $\delta\in\cS\cap\bigcap\limits_{j<\delta}C_j$
(so in particular $\delta\in E_0$) is such that
\begin{enumerate}
\item[$(*)_\delta$]\qquad\qquad $\langle h^{[\zeta]}\circ F_\delta(\alpha):
\alpha<\delta\rangle= \langle r^-_\alpha:\alpha<\delta\rangle$.
\end{enumerate}
We are going to show that $q_\delta\leq r_\delta$ and for this we prove by
induction on $\varepsilon\in(\zeta+1)\cap N$ that $q_\delta\restriction
\varepsilon\leq r_\delta\restriction\varepsilon$. For $\varepsilon=\zeta$
this is the desired conclusion.

For $\varepsilon=0$ this is trivial, and for a limit $\varepsilon$ it
follows from the definition of the order (and the inductive hypothesis).

So assume that we have proved $q_\delta\restriction\varepsilon\leq
r_\delta\restriction\varepsilon$, $\varepsilon<\zeta$, and let us consider 
the restrictions to $\varepsilon+1$. If $\Upsilon(\varepsilon)\geq\delta$
then by the choice of conditions $s,s^-$ in Case 1, we know that 
\[r_\delta\restriction\varepsilon\forces\mbox{`` }(\forall i<\delta)(\exists
j'<\delta)(p_i(\varepsilon)\leq r^-_{j'}(\varepsilon))\mbox{ ''.}\]
Now look at the clause (iii) of the choice of the $p_\delta$ at the
beginning: what we have just stated (and $(*)_\delta$) implies that
\[r_\delta\restriction\varepsilon\forces\mbox{`` }p_\delta(\varepsilon)
\mbox{ is an upper bound to }\{q_\delta(\varepsilon)\}\cup\{p_i(\varepsilon)
:i<\delta\}\mbox{ ''.}\]
thus, $r_\delta\restriction\varepsilon\forces\mbox{`` }q_\delta(\varepsilon)
\leq p_\delta(\varepsilon)\leq r_\delta(\varepsilon)\mbox{ ''}$, so we are
done. Suppose now that $\Upsilon(\varepsilon)<\delta$ and let $j_1=\min(\cS'
\setminus\delta)$. Look at what the generic player has written aside:
$r_\delta\restriction\varepsilon$ forces that $\langle r^\ominus_{j'}
(\varepsilon),r^\oplus_{j'}(\varepsilon),\mathop{\triangle}\limits_{\xi<
\lambda}\name{C}_j^\xi(\varepsilon):j<j_1\rangle$ is a play according to
$\st_\varepsilon$ and $\delta\in\bigcap\limits_{j,\xi<\delta}\name{C}_j^\xi
(\varepsilon)\cap\cS^{\langle\delta\rangle}$, so we are clearly done in this
case too (remember the choice of $r^*$).   
\end{proof}
\end{proof}

\begin{remark}
Note that if the iterands $\name{\bbQ}_\xi$ are (forced to be)
$\lambda$--lub--closed, then the proof of \ref{1.6} substantially
simplifies. 
\end{remark}

\section{Examples}
Our first example of a proper over $\lambda$ forcing notion is a relative of
the forcing introduced by Baumgartner for adding a club to $\aleph_1$. Its
variants were also studied in Abraham and Shelah \cite{AbSh:146}; see also
\cite[Ch.III]{Sh:f}.    

The forcing notion $\bbP^*$ is defined as follows:
\smallskip

\noindent{\bf a condition in $\bbP^*$}\quad is a function $p$ such that 
\begin{enumerate}
\item[(a)] $\Dom(p)\subseteq\lambda^+$, $\Rang(p)\subseteq\lambda^+$,
$|\Dom(p)|<\lambda$, and 
\item[(b)] if $\alpha_1<\alpha_2$ are both from $\Dom(p)$, then
$p(\alpha_1)<\alpha_2$; 
\end{enumerate}
\noindent{\bf the order $\leq$ of $\bbP^*$}\quad is the inclusion
$\subseteq$. 
\smallskip

Clearly,
\begin{proposition}
$\bbP^*$ is $\lambda$--lub--complete and $|\bbP^*|=\lambda^+$.
\end{proposition}

But also,
\begin{proposition}
\label{pstar}
$\bbP^*$ is proper over $\lambda$.
\end{proposition}

\begin{proof}
Assume $N\prec({\mathcal H}(\chi),{\in},{<^*_\chi})$ is as in \ref{pre1.2},
and $p\in\bbP^*\cap N$. 

Put $j^*=N\cap\lambda^+$ and $r=p\cup\{\langle j^*,j^*\rangle\}$.

\begin{claim}
\label{clx1}
\begin{enumerate}
\item If $r'\in\bbP$, $r\leq r'$, then $r'\restriction j^*\in N\cap\bbP^*$
and $r'\restriction j^*\leq r'$.  
\item If $r'\in\bbP$, $r\leq r'$, and $r'\restriction j^*\leq r''\in N\cap
\bbP^*$, then $r'\cup r''\in \bbP^*$ is stronger than both $r'$ and $r''$.
\item If $\bar{p}=\langle p_\xi:\xi<{\zeta^*}\rangle\subseteq\bbP^*$ is
$\leq$--increasing and ${\zeta^*}<\lambda$, then $\bar{p}$ has a least upper
bound $q\in\bbP^*$, and $q\restriction j^*$ is a least upper bound of
$\langle p_\xi\restriction j^*:\xi<{\zeta^*}\rangle$. 
\end{enumerate}
\end{claim}

\begin{proof}[Proof of the Claim]
(1)\quad By the definition of $\bbP^*$, 
\[\alpha\in\Dom(r')\cap j^*\quad\Rightarrow\quad r'(\alpha)<j^*
\quad\Rightarrow\quad r'(\alpha)\in N.\] 
(2), (3)\quad Should be clear.
\end{proof}

\begin{claim}
\label{clx2}
$r$ is $N$--generic for semi--diamonds (see \ref{gencon}).
\end{claim}

\begin{proof}[Proof of the Claim]
Suppose that $D$ is a normal filter on $\lambda$, ${\mathcal S}\in D^+$.
Let $h:\lambda\longrightarrow N$ be such that $N\cap\bbP^*\subseteq
\Rang(h)$, $\bar{F}=\langle F_\delta:\delta\in {\mathcal S}\rangle$ be a
$(D,\cS)$--semi diamond, and let $\bar{q}=\langle q_\delta:\delta\in 
{\mathcal S}\rangle$ be an $(N,h,\bbP^*,\bar{F})$--candidate.

We have to show that the condition $r$ is $(N,h,\bbP^*)$--generic for
$\bar{q}$ over $\bar{F}$, and for this we have to show that the generic
player has a winning strategy in the game $\Game(r,N,h,\bbP^*,\bar{F},
\bar{q})$. Note that the set 
\[E_0\stackrel{\rm def}{=}\{\delta<\lambda:\delta\mbox{ is a limit of
members of }\cS\ \}\]
is a club of $\lambda$ (so $E_0\in D$). Now, the strategy that works for the
generic player is the following one:   
\medskip

\noindent At stage $\delta\in{\mathcal S}$ of the play, when a sequence
$\langle r^-_i,r_i,C_i:i<\delta\rangle$ has been already constructed, the
generic player lets $C_\delta=E_0\setminus(\delta+1)$ and it asks:
\begin{enumerate}
\item[$(*)$] Is there a common upper bound to $\{r_i:i<\delta\}\cup\{
q_\delta\}$ ?
\end{enumerate}
If the answer to $(*)$ is ``yes'', then the generic player puts $Y=\{r_i:i
<\delta\}\cup\{q_\delta\}$; otherwise it lets $Y=\{r_i:i<\delta\}$. Now
it chooses $r_\delta$ to be the $<^*_\chi$--first element of $\bbP^*$
stronger than all members of $Y$ and $r^-_\delta=r_\delta\restriction j^*\in
N$. 
\medskip

Why the strategy described above is the winning one? Let $\langle r_i^-,
r_i,C_i:i<\lambda\rangle$ be a play according to this strategy. Suppose that
$\delta\in\cS\cap\mathop{\triangle}\limits_{i<\lambda} C_i$ is a limit
ordinal such that $\langle h\circ F_\delta(\alpha):\alpha<\delta 
\rangle=\langle r^-_\alpha:\alpha<\delta\rangle$. So, $q_\delta$ is
stronger than all $r^-_\alpha$ (for $\alpha<\delta$), and for cofinally many
$\alpha<\delta$ we have $r^-_\alpha=r_\alpha\restriction j^*$. Therefore,
$q_\delta\geq r^-_\alpha\restriction j^*$ and (by \ref{clx1}) $\{r_\alpha:
\alpha<\delta\}\cup\{q_\delta\}$ has an upper bound. Now look at the choice
of $r_\delta$. 
\end{proof}
The proposition follows immediately from \ref{clx2}.
\end{proof}

\begin{proposition}
\begin{enumerate}
\item $\bbP^*$ is $\alpha$--proper over $\lambda$ if and only if
$\alpha<\lambda$. 
\item If $D$ is a normal filter on $\lambda$, ${\mathcal S}\in D^+$, and
$\bar{F}$ is a $(D,{\mathcal S})$--diamond, $0<\alpha<\lambda^+$, then
$\bbP^*\in K^{\alpha,s}_D[\bar{F}]$ if and only if $\alpha<\lambda$.
\end{enumerate}
\end{proposition}

\begin{proof}
(1)\quad Follows from (2).
\medskip

\noindent (2)\quad Assume $\alpha<\lambda$.

Let $\bar{N}=\langle N_\beta:\beta<\alpha\rangle$, $h_\beta:\lambda
\longrightarrow N_\beta$ and $\bar{q}^\beta$ be as in \ref{1.3}(1b), $p\in
\bbP^*\cap N_0$. Let $j^*_\beta=N_\beta\cap\lambda^+$ (for $\beta<\alpha$)
and put $r=p\cup\{(j^*_\beta,j^*_\beta):\beta<\alpha\}$. Clearly $r\in
\bbP^*$ and $r\restriction j^*_\beta\in N_\beta$ for each $\beta<\alpha$
(remember $\bar{N}\restriction\beta\in N_\beta$). By the proof of
\ref{pstar}, the condition $r\restriction j^*_{\beta+1}$ is $(N_\beta,
h_\beta,\bbP^*)$--generic for $\bar{q}^\beta$ over $\bar{F}$
\medskip

To show that $\bbP^*\notin K^{\alpha,s}_D[\bar{F}]$ for $\alpha\geq\lambda$,
it is enough to do this for $\alpha=\lambda$. So, pick any $\bar{N}=\langle
N_\beta:\beta<\lambda\rangle$, $h_\beta:\lambda\longrightarrow N_\beta$ and
$\bar{q}^\beta$ as in \ref{1.3}(1b), and let $N_\lambda=\bigcup\limits_{
\alpha<\lambda} N_\alpha$. 

Let $\name{\varphi}$ be a $\bbP^*$--name for the generic partial function
from $\lambda^+$ to $\lambda^+$, that is, $\forces_{\bbP^*}\name{\varphi}=
\bigcup\name{G}_{\bbP^*}$. We claim that
\begin{enumerate}
\item[$(\circledast)$] $\forces_{\bbP^*}\mbox{`` }(\exists\beta<\lambda)
(\exists i\in\Dom(\name{\varphi})\cap N_\beta)(\name{\varphi}(i)\notin
N_\beta)\mbox{ ''.}$
\end{enumerate}
Why? Let $p\in\bbP^*$. Take $\beta_0<\lambda$ such that $\Dom(p)\cap
N_\lambda\subseteq N_{\beta_0}$ (remember $|p|<\lambda$). If for some $i\in
\Dom(p)\cap N_{\beta_0}$ we have $p(i)\notin N_{\beta_0}$, then 
\[p\forces\mbox{`` }(\exists i\in\Dom(\name{\varphi})\cap N_{\beta_0})(
\name{\varphi}(i)\notin N_{\beta_0})\mbox{ ''.}\]
Otherwise, we let $\delta^*=N_{\beta_0}\cap\lambda^+$ and $\delta^{**}=
N_{\beta_0+1}\cap\lambda^+$, and we put $q=p\cup\{(\delta^*,\delta^{**})
\}$. Then clearly $q\in\bbP^*$ is a condition stronger than $p$ and
\[q\forces\mbox{`` }(\exists i\in\Dom(\name{\varphi})\cap N_{\beta_0+1})(
\name{\varphi}(i)\notin N_{\beta_0+1}).\]

It should be clear that $(\circledast)$ implies that there is no condition
$r\in\bbP^*$ which is $(N_\beta,h_\beta,\bbP^*)$--generic for
$\bar{q}^\beta$ for all $\beta<\alpha$ (remember \ref{1.2A}(2)).
\end{proof}

For the second example we assume the following.
\begin{context}
\label{conlast}
\begin{enumerate}
\item[(a)] $\lambda,D,\cS,\cS'$ are as in \ref{1.1},
\item[(b)] $\cS^*\subseteq\cS^{\lambda^+}_\lambda\stackrel{\rm def}{=}\{
\delta<\lambda^+:\cf(\delta)=\lambda\}$, 
\item[(c)] $\langle A_\delta,h_\delta:\delta\in\cS^*\rangle$ is such that
for each $\delta\in\cS^*$:
\item[(d)] $A_\delta\subseteq\delta$, $\otp(A_\delta)=\lambda$ and
$A_\delta$ is a club of $\delta$, and 
\item[(e)] $h_\delta:A_\delta\longrightarrow\lambda$.
\end{enumerate}
\end{context}
The forcing notion $\bbQ^*$ is defined as follows:
\smallskip

\noindent{\bf a condition in $\bbQ^*$}\quad is a tuple $p=(u^p,v^p,
\bar{e}^p,h^p)$ such that 
\begin{enumerate}
\item[(a)] $u^p\in [\lambda^+]^{<\lambda}$, $v^p\in [\cS^*]^{<\lambda}$, 
\item[(b)] $\bar{e}^p=\langle e^p_\delta:\delta\in v^p\rangle$, where each
$e^p_\delta$ is a closed bounded subset of $A_\delta$, and $e^p_\delta
\subseteq u^p$, and
\item[(c)] if $\delta_1<\delta_2$ are from $v^p$, then 
\[\sup(e_{\delta_2})>\delta_1\quad\mbox{ and }\quad \sup(e_{\delta_1})>
\sup(A_{\delta_2}\cap\delta_1),\] 
\item[(d)] $h^p:u^p\longrightarrow\lambda$ is such that for each $\delta\in
v^p$ we have
\[h^p\restriction \{\alpha\in e_\delta:\otp(\alpha\cap e_\delta)\in\cS'\}
\subseteq h_\delta;\]
\end{enumerate}
\noindent{\bf the order $\leq$ of $\bbQ^*$}\quad is such that $p\leq q$ if
and only if $u^p\subseteq u^q$, $h^p\subseteq h^q$, $v^p\subseteq v^q$, and
for each $\delta\in v^p$ the set $e^q_\delta$ is an end-extension of
$e^p_\delta$.
\smallskip

A tuple $p=(u^p,v^p,\bar{e}^p,h^p)$ satisfying clauses (a), (b) and (d)
above will be called {\em a pre-condition}. Note that every pre-condition
can be extended to a condition in $\bbQ^*$.
\smallskip

Plainly:
\begin{proposition}
The forcing notion $\bbQ^*$ is $\lambda$--lub--complete. Also $\bbQ^*$
satisfies the $\lambda^+$--chain condition. 
\end{proposition}

\begin{proposition}
$\bbQ^*$ is proper over $\lambda$.
\end{proposition}

\begin{proof}
Assume $N\prec({\mathcal H}(\chi),{\in},{<^*_\chi})$ is as in \ref{pre1.2},
$\langle A_\delta,h_\delta:\delta\in\cS^*\rangle\in N$ and $p\in\bbQ^*\cap
N$. We are going to show that the condition $p$ is $N$--generic for
semi-diamonds. 

So suppose that $h,\bar{q}$ and $\bar{F}$ are as in \ref{pre1.2}. For
$r\in\bbQ^*$, let $r\restriction N$ be such that $u^{r\restriction N}=u \cap
N$, $v^{r\restriction N}=v\cap N$, $\bar{e}^{r\restriction N}= \bar{e}^r
\restriction N$, $h^{r\restriction N}=h^r\restriction N$. Note that
$r\restriction N\in N$. Let us describe the winning strategy of the generic
player in the game $\Game(p,N,h,\bbP^*,\bar{F},\bar{q})$. For this we first
fix a list $\{j_i:i<\lambda\}$ of $N\cap \cS^*$, and we let
$E_0=\{\delta<\lambda:\delta$ is a limit of members of $\cS\ \}$.

Suppose that we arrive to a stage $\delta\in\cS$ and $\langle r^-_i,
r_i,C_i:i<\delta\rangle$ is the sequence played so far. The generic player
first picks a condition $r_\delta'$ stronger than all $r_i$'s played so far
and, if possible, stronger than $q_\delta$. Then it plays a condition
$r_\delta$ above $r_\delta'$ such that
\begin{itemize}
\item if $\gamma\in v^{r_\delta}$, then $\otp(e^{r_\delta}_\gamma)>\delta$,
and  
\item $\{j_i:i<\delta\}\subseteq v^{r_\delta}$,
\end{itemize}
and $r^-_\delta=r_\delta\restriction N$. The set $C_\delta$ played a this
stage is $[\alpha,\lambda)\cap E_0$, where $\alpha$ is the first ordinal
such that $v^{r_\delta}\cap N\subseteq\{j_i:i<\alpha 
\}$ and $\otp(A_\gamma\cap (\max(e^{r_\delta}_\gamma)+1))<\alpha$ for all
$\gamma\in v^{r_\delta}$. 

Why is this a winning strategy? Let $\langle r^-_i,r_i,C_i:i<\lambda\rangle$
be a play according to this strategy, and suppose that $\delta\in\cS\cap
\mathop{\triangle}\limits_{i<\lambda} C_i$ is a limit ordinal such that 
\[\langle h\circ F_\delta(\alpha):\alpha<\delta\rangle=\langle
r^-_\alpha:\alpha<\delta\rangle.\]
Note that then
\begin{enumerate}
\item[(i)] if $\gamma\in\bigcup\limits_{i<\delta} v^{r_i}$ then
$\bigcup\limits_{i<\delta} e^{r_i}_\gamma$ is an unbounded subset of
$\{\varepsilon\in A_\gamma:\otp(\varepsilon\cap A_\gamma)<\delta\}$, and 
\item[(ii)] $\bigcup\limits_{i<\delta}v^{r_i}\cap N=\{j_i:i<\delta\}$.
\end{enumerate}
We want to show that there there is a common upper bound to $\{r_i:i<\delta
\}\cup\{q_\delta\}$ (which, by the definition of our strategy, will finish
the proof). First we choose a pre-condition $r=(u^r,v^r,\bar{e}^r,h^r)$ such
that:
\begin{itemize}
\item $v^r=v^{q_\delta}\cup\bigcup\limits_{i<\delta} v^{r_i}$,
\item if $\gamma\in v^{q_\delta}$, then we let $e^r_\gamma=
e^{q_\delta}_\gamma$, if $\gamma\in\bigcup\limits_{i<\delta}v^{r_i}\setminus
v^{q_\delta}$, then 
\[e^r_\gamma=\bigcup\{e^{r_i}_\gamma:i<\delta,\ \gamma\in v^{r_i}\}\cup\{
\mbox{ the $\delta^{\rm th}$ member of $A_\gamma$ }\},\]
\item $u^r=u^{q_\delta}\cup\bigcup\limits_{i<\delta} u^{r_i}\cup
\{\mbox{the $\delta^{\rm th}$ member of $A_\gamma$}:\gamma\in v^r\setminus
v^{q_\delta}\}$,
\item $h^r\supseteq h^{q_\delta}\cup\bigcup\limits_{i<\delta} h^{r_i}$.
\end{itemize}
Why is the choice possible? As $\delta\notin\cS'$ ! Now we may extend $r$ to
a condition in $\bbQ^*$ picking for each $\gamma\in v^r$ large enough
$\beta_\gamma\in A_\gamma$ and adding $\beta_\gamma$ to $e^r_\gamma$ (and
extending $u^r,h^r$ suitably). 
\end{proof}
\medskip

Our last example is a natural generalization of the forcing notion
$\bbD_\omega$ from Newelski and Ros{\l}anowski \cite{NeRo93}. Let us work in
the context of \ref{1.1}.

\begin{definition}
\begin{enumerate}
\item A set $T\subseteq{}^{<\lambda}\lambda$ is {\em a complete
$\lambda$--tree\/} if
\begin{enumerate}
\item[(a)] $(\forall\eta\in T)(\exists\nu\in T)(\eta\vartriangleleft\nu)$,
and $T$ has the $\vartriangleleft$--smallest element called $\mrot(T)$, 
\item[(b)] $(\forall\eta,\nu\in {}^{<\lambda}\lambda)(\mrot(T)
\trianglelefteq \eta\vartriangleleft\nu\in T\ \Rightarrow\ \eta\in T)$,
\item[(c)] if $\langle\eta_i:i<\delta\rangle\subseteq T$ is a
$\vartriangleleft$--increasing chain, $\delta<\lambda$, then there is
$\eta\in T$ such that $\eta_i\vartriangleleft\eta$ for all $i<\delta$. 
\end{enumerate}
Let $T\subseteq {}^{<\lambda}\lambda$ be a complete $\lambda$--tree.
\item For $\eta\in T$ we let $\suc_T(\eta)=\{\alpha<\lambda: \eta^\frown
\langle\alpha\rangle\in T\}$.
\item We let $\splt(T)=\{\eta\in T:|\suc_T(\eta)|>1\}$.
\item A sequence $\rho\in{}^\lambda\lambda$ is a {\em $\lambda$--branch
through $T$} if 
\[(\forall\alpha<\lambda)({\rm lh}(\mrot(T))\leq\alpha\quad \Rightarrow
\quad\rho\restriction\alpha\in T).\]
The set of all $\lambda$--branches through $T$ is called $\lim_\lambda(T)$.  
\item A subset $\bF$ of the $\lambda$--tree $T$ is {\em a front} in $T$ if no
two distinct members of $\bF$ are $\vartriangleleft$--comparable and 
\[(\forall\eta\in{\lim}_\lambda(T))(\exists\alpha<\lambda)(\eta\restriction 
\alpha\in\bF).\] 
\item For $\eta\in T$ we let $(T)^{[\eta]}=\{\nu\in T:\eta\vartriangleleft
\nu\}$. 
\end{enumerate}
\end{definition}

Now we define a forcing notion $\bbD_\lambda$:\\
{\bf A condition in $\bbD_\lambda$}\quad is a complete $\lambda$--tree $T$
such that
\begin{enumerate}
\item[(a)] $\mrot(T)\in\splt(T)$ and $(\forall\eta\in \splt(T))(\suc_T(
\eta)=\lambda)$,
\item[(b)] $(\forall\eta\in T)(\exists\nu\in T)(\eta\vartriangleleft\nu\in
\splt(T))$, 
\item[(c)] if $\delta<\lambda$ is limit and a sequence $\langle\eta_i:i<
\delta\rangle\subseteq\splt(T)$ is $\vartriangleleft$--increasing,
\underline{then} $\eta=\bigcup\limits_{i<\delta}\eta_i\in\splt(T)$.
\end{enumerate}
{\bf The order of $\bbD_\lambda$}\quad is the reverse inclusion.

\begin{proposition}
$\bbD_\lambda$ is proper over $\lambda$.
\end{proposition}

\begin{proof}
First let us argue that $\bbD_\lambda$ is $\lambda$--lub--complete. So
suppose that $T_\alpha\in \bbD_\lambda$ are such that $(\forall\alpha<
\beta<\delta)(T_\beta\subseteq T_\alpha)$, $\delta<\lambda$. We claim that
$T\stackrel{\rm def}{=}\bigcap\limits_{\alpha<\delta}T_\alpha$ is a
condition in $\bbD_\lambda$. Clearly $T$ is a tree, and $\mrot(T)=
\bigcup\limits_{\alpha<\delta}\mrot(T_\alpha)$. By clause (c) (for
$T_\alpha$'s) we see that $\suc_T(\mrot(T))=\lambda$, and in a similar way
we justify that $T$ satisfies other demands as well. 

Now suppose that $D,\cS,N,h,\bar{F}$ and $\bar{q}$ are as in \ref{pre1.2},
$T\in\bbD_\lambda\cap N$. Choose inductively complete $\lambda$--trees
$T_\alpha\in \bbD_\lambda\cap N$ and fronts $\bF_\alpha\subseteq T_\alpha$
(of $T_\alpha$) such that
\begin{enumerate}
\item[(i)]   $\mrot(T_\alpha)=\mrot(T)$,
\item[(ii)]  if $\alpha\leq\beta<\lambda$, then $T_\beta\subseteq
T_\alpha\subseteq T$ and $\bF_\alpha\subseteq \splt(T_\beta)$, and 
\item[(iii)] if $\eta\in\bF_\alpha$ then $\otp(\{i<{\rm lh}(\eta):
\eta\restriction i\in\splt(T_\alpha)\})=\alpha$,
\item[(iv)]  if $\delta$ is limit, then $T_\delta=\bigcap\limits_{\alpha
<\delta}T_\alpha$, 
\item[(v)]   if $\delta\in\cS$ is limit and $\langle h\circ
F_\delta(\alpha):\alpha<\delta\rangle$ is an increasing sequence of
conditions from $\bbD_\lambda\cap N$ and $\bigcap\limits_{\alpha<\delta}
h\circ F_\delta(\alpha)\subseteq T_\delta$, and
$\eta=\bigcup\limits_{\alpha<\delta}\mrot(h\circ F_\delta(\alpha))\in
\bF_\delta$, \underline{then} for some $\nu\in T_\delta$ we have
\[\eta\vartriangleleft\nu\in\bF_{\delta+1}\quad\mbox{ and }\quad
q_\delta\leq (T_{\delta+1})^{[\nu]}.\]
\end{enumerate}
Now we let $r=\bigcap\limits_{\alpha<\lambda} T_\alpha$. It should be
clear that $r\in\bbD_\lambda$.

\begin{claim}
$r$ is $(N,h,\bbD_\lambda)$--generic for $\bar{q}$ over $\bar{F}$.
\end{claim}

\begin{proof}[Proof of the Claim]
We have to describe a winning strategy of the generic player in the game 
$\Game(r,N,h,\bbD_\lambda,\bar{F},\bar{q})$. Let $E_0$ be the club of limits
of members of $\cS'$. Let the generic player play as follows.

Assume we have arrived to stage $i\in\cS$ of the play when $\langle
r^-_j,r_j,C_j:j<i\rangle$ has been already constructed. If $i\notin E_0$
then the generic player chooses $r^-_i, r_i\in\bbD_\lambda$ such that
\begin{enumerate}
\item[(A)] $r_i\subseteq\bigcap\limits_{j<i}r_j$, $r_i^-\subseteq
\bigcap\limits_{j<i}r_j^-\cap \bigcap\limits_{j<i} T_j$, and $r^-_i\in N$,
$r^-_i\leq r_i$, 
\item[(B)] $\mrot(r^-_i)=\mrot(r_i)\in \bF_{\alpha(i)}$ for some $\alpha(i)
>i$,
\end{enumerate}
and lets $C_i=E_0\setminus (\alpha(i)+1)$. If $i\in E_0$ then the generic
player picks $r_i,r^-_i$ satisfying (A) + (B) and such that 
\begin{enumerate}
\item[(C)] if possible, then $q_\delta\leq r^-_i$
\end{enumerate}
and it takes $C_i$ as earlier.

Why is this a winning strategy? First, as $\bbD_\lambda$ is
$\lambda$--lub--complete, the play really lasts $\lambda$ moves. Suppose that
$\delta\in \cS\cap\bigcap\limits_{i<\delta}C_i$ is such that 
\[\langle h\circ F_\delta(\alpha):\alpha<\delta\rangle =\langle r^-_\alpha:
\alpha<\delta\rangle.\]
Let $\eta=\bigcup\limits_{\alpha<\delta}\mrot(r^-_\alpha)$. Note that (as
$\delta\in E_0$ and by (B)) we have $\eta\in\bF_\delta$ and (by (A))
$\bigcap\limits_{\alpha<\delta}r^-_\alpha$ is included in
$T_\delta$. Therefore, by clause (v) of the choice of the $T_\delta$, 
for some $\nu\in T_\delta$ we have $\eta\vartriangleleft\nu\in
\bF_{\delta+1}$ and $q_\delta\leq (T_{\delta+1})^{[\nu]}$. But this
immediately implies that it was possible to choose $r^-_i$ stronger than
$q_\delta$ in (C) (remember $r=\bigcap\limits_{\alpha<\delta} T_\alpha$). 
\end{proof}

\end{proof}

\section{Discussion}

\subsection{The Axiom}\label{axiom} We can derive Forcing Axiom as usual,
see \cite[Ch. VII, VIII]{Sh:f}. E.g., if $\kappa$ is a supercompact cardinal
larger than $\lambda$, then we can find a $\kappa$--cc $\lambda$--complete,
proper over $D$--semi diamonds forcing notion $\bbP$ of cardinality $\kappa$
such that 
\begin{itemize}
\item $\forces_{\bbP}2^\lambda=\kappa$
\item $\bbP$ collapses every $\mu\in (\lambda^+,\kappa)$, no other cardinal
is collapsed,
\item in $\bV^\bbP$:\\
\underline{if} $\bbQ$ is a forcing notion proper over $D$--semi diamonds,
$\cI_\alpha$ are open dense subsets of $\bbQ$ for $\alpha<\lambda^+$,\\
\underline{then} there is a directed set $G\subseteq\bbQ$ intersecting every
$\cI_\alpha$ (for $\alpha<\lambda^+$).
\end{itemize}
If we restrict ourselves to $|\bbQ|=\kappa$, it is enough that $\kappa$ is
indiscrebable enough. 

In ZFC, we have to be more careful concerning $\bbQ$.

\subsection{Future applications} Real applications of the technology
developed here will be given in a forthcoming paper Ros{\l}anowski and
Shelah \cite{RoSh:777}, where we will present more examples of proper for
$\lambda$ forcing notions (concentrating on the case of inaccessible
$\lambda$). We start there developing a theory parallel to that of
Ros{\l}anowski and Shelah \cite{RoSh:470}, \cite{RoSh:628}, \cite{RoSh:672}
aiming at generalizing many of the cardinal characteristics of the continuum
to larger cardinals.  

\subsection{Why our definitions?} The main reason why our definitions are
(perhaps) somewhat complicated is that, in addition to ZFC limitations, we
wanted to cover some examples with ``large creatures'' (to be presented in
\cite{RoSh:777}). We also wanted to have a real preservation theorem: the
(limit of the) iteration is of the same type as the iterands (though for
many applications the existence of $(N,\bbP_\zeta)$--generic conditions
could be enough). 

Why do we have the sets $C_i$ in the game, and not just say that ``the set
of good $\delta$'s is in $D$''? It is caused by the fact that already if we
want to deal with the composition of two forcing notions (the successor
step), the respective set from $D$ would have appeared only {\em after} the
play, and there would be simply too many possible sets to consider. With the 
current definition the generic player discovers during the play which
$\delta\in\cS$ are active. 

Why semi--diamonds (and not just diamonds)? As we want that
$\bar{q}^{\langle i\rangle},\bar{F}^{\langle i\rangle}$ are as claimed in
\ref{cl2} (for the respective parameters). 

\subsection{Strategic completeness} We may replace ``$\lambda$--complete''
by (a variant of) ``strategically $\lambda$--complete''. This requires some
changes in our definitions (and proofs) and it will be treated in
\cite{Sh:F509}. 

\subsection{Relation to \cite{Sh:587}}

There is a drawback in the approach presented in this paper: we do not
include the one from \cite{Sh:587}, say when ${\mathcal S}\subseteq
{\mathcal S}^{\lambda^+}_\lambda$ is stationary and ${\mathcal
S}^{\lambda^+}_\lambda\setminus {\mathcal S}$ is also stationary.

One of possible modification of the present definitions for the case
of inaccessible $\lambda$, can be sketched as follows. We have
$\langle\lambda_\delta:\delta\in{\mathcal S}\rangle$, $\lambda_\delta=
(\lambda_\delta)^{|\delta|}$; $\bar{q}=\langle q_\delta:\delta\in{\mathcal
S} \rangle$  is replaced by $\bar{q}=\langle q_{\delta,t}:\delta\in{\mathcal
S},\ t\in{\rm Par}_{\delta^*,\delta}\rangle$ (where $\delta^*=N\cap
\lambda^+$), and $\bar{\rm Par}=\langle{\rm Par}_{\delta^*,\delta}:\delta\in
{\mathcal S}\rangle\in\bV$ is constant for the iteration (like $D$). 

In the forcing $\bbP$: for $\bar{p}=\langle p_j:j<\delta\rangle$,
$\delta\in{\mathcal S}$, $t\in {\rm Par}_{\delta^*,\delta}$, there is an
upper bound $q[\bar{p},t]$ of $\bar{p}$ (this is a part of $\bbP$).

For each $\delta$, each ${\rm Par}_{\delta^*,N_\delta\cap\gamma}$,
$\prod\limits_{i\in N_\delta}{\rm Par}_{\delta^*,\delta}$ has cardinality
$\lambda_\delta=(\lambda_\delta)^{|\delta|}$ ($N_\delta$ is of cardinality
$|\delta|$; $\gamma$ is the length of the iteration). Having $\langle
p_j:j<\delta\rangle\subseteq N_\delta$ we can find $\langle q^\delta_t: t\in
{\rm Par}_{\delta,N_\delta\cap\gamma}\rangle$ as in \cite{Sh:587}. 

Several (more complex) variants of properness over semi--diamonds will be
presented in \cite{Sh:F509} and Ros{\l}anowski and Shelah \cite{RoSh:777}.


\end{document}